\documentclass[12pt]{article}

\usepackage{amsthm,amsmath,amsfonts,epsfig,amssymb,latexsym}
\usepackage{mathrsfs} 
\usepackage[T1]{fontenc}
\usepackage{graphicx}
\usepackage{enumerate}
\usepackage{xy,xypic}
\usepackage{graphics}
\usepackage{footnote}
\usepackage{scalefnt}
\usepackage{tikz}
\usepackage{color}
\usetikzlibrary{shapes,arrows,fit,calc,positioning}
\usetikzlibrary{matrix}
\xyoption{all}
\title{Witt, $GW$, $K$-theory of quasi-projective schemes}

\author{
 Satya Mandal\footnote{Partially supported by a General Research Grant (no 2301857) from U. of Kansas}
 \\ 
{\small University of Kansas, Lawrence, Kansas 66045;}
{\small {\it  mandal@ku.edu 
  }}\\
 }

\begin{document}
\renewcommand{\baselinestretch}{1.255}
\setlength{\parskip}{1ex plus0.5ex}
\date{9 September 2015}
\newcommand{\iso}{\stackrel{\sim}{\longrightarrow}}
\newcommand{\sur}{\twoheadrightarrow}

\newcommand{\eop}{\hfill \rule{2mm}{2mm}}
\newcommand{\pf}{\noindent{\bf Proof.~}}
\newcommand{\PD}{\dim_{{\SA}}}
\newcommand{\PDV}{\dim_{{\SV(X)}}}

\newcommand{\lra}{\longrightarrow}
\newcommand{\hra}{\hookrightarrow} 
\newcommand{\Lra}{\Longrightarrow}
\newcommand{\Llra}{\Longleftrightarrow}
\newcommand{\pic}{The proof is complete.~}
\newcommand{\Sp}{\mathrm{Sp}}

\newcommand{\Dia}{\diagram}

\newcommand{\bE}{\begin{enumerate}}
\newcommand{\eE}{\end{enumerate}}

\newtheorem{theorem}{Theorem}[section]
\newtheorem{proposition}[theorem]{Proposition}
\newtheorem{lemma}[theorem]{Lemma}
\newtheorem{definition}[theorem]{Definition}
\newtheorem{corollary}[theorem]{Corollary}
\newtheorem{construction}[theorem]{Construction}
\newtheorem{notations}[theorem]{Notations}
\newtheorem{remark}[theorem]{Remark}
\newtheorem{question}[theorem]{Question}
\newtheorem{example}[theorem]{Example} 

\newcommand{\bD}{\begin{definition}}
\newcommand{\eD}{\end{definition}}
\newcommand{\bP}{\begin{proposition}}
\newcommand{\eP}{\end{proposition}}
\newcommand{\bL}{\begin{lemma}}
\newcommand{\eL}{\end{lemma}}
\newcommand{\bT}{\begin{theorem}}
\newcommand{\eT}{\end{theorem}}
\newcommand{\bC}{\begin{corollary}}
\newcommand{\eC}{\end{corollary}} 
\newcommand{\TCP}{\textcolor{purple}}
\newcommand{\TCM}{\textcolor{magenta}}
\newcommand{\TCR}{\textcolor{red}}
\newcommand{\TCB}{\textcolor{blue}}
\newcommand{\TCG}{\textcolor{green}}
\def\spec#1{\mathrm{Spec}\left(#1\right)}
\def\proj#1{\mathrm{Proj}\left(#1\right)}

\def\m{\mathfrak {m}}
\def\CA{\mathcal {A}}
\def\CB{\mathcal {B}}
\def\CP{\mathcal {P}}
\def\CC{\mathcal {C}}
\def\CD{\mathcal {D}}
\def\CE{\mathcal {E}}
\def\CF{\mathcal {F}}
\def\CG{\mathcal {G}}
\def\CH{\mathcal{H}}
\def\CK{\mathcal{K}}
\def\CL{\mathcal{L}}
\def\CM{\mathcal{M}}
\def\CN{\mathcal{N}}
\def\CO{\mathcal{O}}
\def\CP{\mathcal{P}}
\def\CQ{\mathcal{Q}}

\def\CR{\mathcal{R}}

\def\CS{\mathcal{S}}
\def\CT{\mathcal{T}}
\def\CU{\mathcal{U}}
\def\CV{\mathcal{V}}
\def\CW{\mathcal{W}} 
\def\CX{\mathcal{X}}
\def\CY{\mathcal{Y}}
\def\CZ{\mathcal{Z}}

\newcommand{\smallcirc}[1]{\scalebox{#1}{$\circ$}}
\def\Z{\mathbb {Z}}
\def\A{\mathbb{A}}
\def\B{\mathbb{B}}
\def\P{\mathbb{P}}
\def\C{\mathbb{C}}
\def\D{\mathbb{D}}
\def\E{\mathbb{E}}
\def\F{\mathbb{F}}
\def\G{\mathbb{G}}

\def\H{\mathbb{H}}
\def\K{\mathbb{K}}
\def\L{\mathbb{L}}
\def\M{\mathbb{M}}
\def\N{\mathbb{N}} 
\def\O{\mathbb{O}}
\def\P{\mathbb{P}} 
\def\R{\mathbb{R}}

\def\SA{\mathscr {A}} 
\def\SB{\mathscr {B}}
\def\SC{\mathscr {C}}
\def\SD{\mathscr {D}}
\def\SE{\mathscr {E}}
\def\SF{\mathscr {F}}
\def\SG{\mathscr {G}}
\def\SH{\mathscr {H}}
\def\SI{\mathscr {I}}
\def\SK{\mathscr {K}} 
\def\SL{\mathscr {L}} 
\def\SM{\mathscr {M}}
\def\SN{\mathscr {N}}
\def\SO{\mathscr {O}}
\def\SP{\mathscr {P}}
\def\SQ{\mathscr {Q}}
\def\SR{\mathscr {R}}
\def\SS{\mathscr {S}}
\def\ST{\mathscr {T}}

\def\SU{\mathscr {U}}

\def\SV{\mathscr {V}}
\def\SX{\mathscr {X}}
\def\SY{\mathscr {Y}}

\def\SZ{\mathscr {Z}} 
 
\maketitle

\vspace{5mm}
\noindent{\bf Abstract}: In this article we continue our investigation 
 of the Derived Equivalences over noetherian quasi-projective schemes $X$,
 over affine schemes $\spec{A}$. 
 For integers $k\geq 0$, let $C\M^k(X)$ denote the category of coherent ${\CO}_X$-modules $\CF$, with
locally free dimension $\PDV(\CF)=k=grade(\CF)$. We prove that there is a zig-zag equivalence 
${\CD}^b\left(C\M^k(X)\right) \to {\SD}^k\left(\SV(X)\right)$ of the derived categories. It follows that there is a sequence 
of zig-zag maps 
$
\diagram 
\K\left(C\M^{k+1}(X)\right) \ar[r] & \K\left(C\M^{k}(X)\right) \ar[r] & \coprod_{x\in X^{(k)}} \K\left(C\M^{k}(X_x)\right) \\
\enddiagram
$
of the $\K$-theory spectra that is a  homotopy fibration.  
In fact, this is 
analogous to the fibrations of the $G$-theory spaces of Quillen (see proof of \cite[Theorem 5.4]{Q}).
We also establish similar homotopy fibrations of ${\bf GW}$-spectra and ${\G}W$-bispectra.

\section{Introduction}

In \cite{Q}, Qullen established the foundation of $K$-theory of regular schemes $X$ in a complete manner. 
In fact, for any scheme $X$, Quillen provides a complete foundation of $K$-theory of 
the category $Coh(X)$ of the coherent sheaves on $X$, along with the filtration of $Coh(X)$  by co-dimension of support of the 
objects $\CF\in Coh(X)$. This   relates to Gersten complexes and spectral sequences associated to any such scheme $X$ (see \cite[\S 5]{Q}). 
The $K$-theory of $Coh(X)$ is also known as $G$-theory.
For regular schemes $X$
 the $K$-theory of the category $\SV(X)$ of locally free sheaves reconciles fully with that of $Coh(X)$. 
 Consequently, the $K$-theory of regular schemes appears very complete.  
 However,  the $K$-theory of non-regular schemes never reached the completeness and harmony 
 that the $K$-theory of regular schemes had achieved. Work of Waldhausen (\cite{W}) and 
Thomason-Trobaugh (\cite{TT}) would be milestones in this respect, most notably for their introduction of derived invariance theorems and localization theorems,
applicable to non-regular schemes. 
Further, while developments in Grothendieck-Witt theory ($GW$-theory)
and Witt theory followed  the foot prints of $K$-theory (\cite{S3, TWG1}), due to lack of any natural duality on $Coh(X)$, the 
situation in these two areas appear even less complete. 
When, $X$ is non-regular, the category $\M(X)$ of coherent sheaves with finite $\SV(X)$-dimension differs from $Coh(X)$. There appears to be a gap in the literature of  $K$-theory, $GW$-theory, and Witt theory, with respect to the place of the category $\M(X)$.  One can speculate, whether this lack of completeness
  is attributable to this gap? The goal of this one and the related articles is to work on this gap and attempt to establish the said literature on non-regular schemes at the same pedestal as that of regular schemes. 
For quasi-projective schemes over noetherian affine schemes, this goal is accomplished to up some degree of satisfaction. The special place of the 
full subcategory $C\M^k(X)\subset \M(X)$ 
  would also be clear subsequently, 
  where for integers $k\geq 0$, $C\M^k(X)$ will denote the full subcategory of objects $\CF$ in $\M(X)$, with $\PDV(\CF)=grade\left(\CF\right)=k$. 



With respect to a certain facets of these three areas, namely Algebraic ${K}$-theory, Grothendieck-Witt ($GW$) theory and Witt theory, a common thread between them is 
their invariance properties with respect equivalences of the associated Derived categories. We review some of the results on such invariances. 
For example, recall the theorem of Thomanson-Trobaugh (\cite[Theorem 1.9.8]{TT}):
suppose ${\bf A}\to {\bf B}$ is a functor of complicial exact categories with weak equivalences. 
Assume that the associated functor
 of the triangulated categories $\ST{\bf A}\to  \ST{\bf B}$ is an equivalence. Then, the induced map ${\bf K}({\bf A}) \to
 {\bf K}({\bf B})$ of the ${\bf K}$-theory spaces is a homotopy equivalence (see \cite[3.2.24]{S1}). The non-connective version of this theorem 
 was given by Schlichting (\cite[Theorem 9]{S2}, also see \cite[3.2.29]{S1}): which states that under the relaxed hypothesis that,
 if $\ST{\bf A}\to   \ST{\bf B}$ is an equivalence up
 to factors, then it induces a homotopy equivalence $\K({\bf A}) \to \K({\bf B})$ of the $\K$-theory spectra. 
 While ${K}$-theory is defined for complicial exact categories with weak equivalences, Schlichting defined 
 Grotherndieck Witt ($GW$)
  spectra and bispectra (\cite{S3}, also see \S \ref{secSpGWKGW})
 of pointed dg categories with weak equivalences and dualities. Invariance theorems of ${\bf GW}$-spectra
 and ${\G}W$-bispectra, similar to $K$-theory,
 were established in \cite[Theorems 6.5, 8.9]{S3}. Contrary to $K$-theory and  $GW$-theory, Balmer defined Witt theory ($W$-theory)  for 
 Triangulated categories with dualities (\cite{TWG1}), which encompasses the Derived categories with dualities. 
 Therefore, the shifted Witt groups   
 are invariant with respect to equivalences of derived categories, by \cite[Theorem 6.2]{TWG1}. Another common thread between these three 
 areas is  the exactness properties of the associated triangulated categories. 
 %
  %
 In particular, the renowned Gerstner complexes, in $K$-theory, $GW$-theory and Witt theory, 
 are obtained by routine manipulation (see Remark \ref{remGerSPSeq}) of the respective invariants, by such derived equivalences and exactness 
 properties of the associated triangulated categories. 
 Therefore, in this article we consider  equivalences of certain derived 
 categories, over quasi-projective schemes, which we state subsequently. 

The readers are referred to (\ref{nota}) for   clarifications regarding notations and 
the definition of grade. 
Other than the notations explained above, 
for integers $k\geq 0$,
$\M^k(X)$ will denote the category
of coherent ${\CO}_X$-modules $\CF$ with finite locally free dimension, and $grade\left(\CF\right) \geq k$.
We prove that, for a noetherian quasi-projective  scheme $X$ over an affine scheme $\spec{A}$, and integers $k\geq 0$,
the functor of the 
 derived categories
$$
\zeta:{\CD}^b\left(C\M^k(X)\right) \to {\CD}^b\left(\M^k(X)\right)
\quad {\rm is~an~equivalence}
$$  
 (see Theorem \ref{mainThmFoxII}).   
We also prove that  the functor of the derived categories
$$
\beta:{\CD}^b\left(\M^{k+1}(X)\right) \to {\CD}^b\left(\M^k(X)\right)
\quad {\rm is~faithfully~full}
$$
(see Theorem \ref{mkp1toMkFlF}). Consequently, the 
functor ${\CD}^b\left(C\M^{k+1}(X)\right) \to {\CD}^b\left(\M^k(X)\right)$ is faithfully full. 
Combining the results in \cite{FM}, we have the following summary of results. 
Consider the commutative diagram 
\begin{equation}\label{f5FunctIntroGWadd}
\diagram
{\CD}^b(C\M^{k+1}(X))\ar[r]^{\zeta}_{\sim}\ar[d]^{\alpha} &{\CD}^b({\M}^{k+1}(X))\ar[r]^{\iota}_{\sim}\ar[d]^{\beta}
& {\SD}^{k+1}(({\M}(X)) \ar[d]^{\gamma}
& {\SD}^{k+1}(\SV(X))\ar[d]^{\eta}\ar[l]_{\quad\iota'}^{\quad\sim}\\ 
{\CD}^b(C\M^k(X))\ar[r]_{\zeta}^{\sim} &{\CD}^b({\M}^k(X))\ar[r]_{\iota}^{\sim}
& {\SD}^k(({\M}(X)) 
& {\SD}^k(\SV(X))\ar[l]^{\quad\iota'}_{\quad\sim}\\\ 
\enddiagram 
\end{equation} 
of functors of derived categories. 
Then, all the horizontal functors are equivalences and, all the  vertical functors are fully faithful (see Theorem \ref{equivThm}). 

Having stated the  equivalence theorem (\ref{equivThm}), we first turn our attention to its consequences to Algebraic 
$K$-theory of quasi-projective schemes $X$ over an affine scheme $\spec{A}$. Note that $C\M^k(X)$ is an exact category. Quillen \cite{Q} defined 
${\bf K}$-theory space ${\bf K}\left(\SE\right)$ of any exact category $\SE$. To incorporate negative $K$-groups, following Bass, Karoubi and others,
 Schlichting formally introduced  (\cite{S1, S2}) $\K$-theory spectrum $\K\left(\SE\right)$ for such exact categories $\SE$. 
 By agreement theorems (\cite[Theorem 1.11.17]{TT}, \cite[3.2.30]{S1}), there   are homotopy equivalences
 $$
 \left\{
 \begin{array}{ll}
  {\bf K}\left(\SE\right) \iso {\bf K}\left(Ch^b\left(\SE\right)\right) &  {\rm of ~the}~{\bf K}{\rm-theory~spaces},\\
 {\K}\left(\SE\right) \iso \K\left(Ch^b\left(\SE\right)\right) &  {\rm of ~the}~\K{\rm-theory~spectra},\\
 \end{array}
 \right.
 $$
 where the right hand sides  correspond to the 
 $K$-theory 
 space/spectrum of the category 
 $Ch^b\left(\SE\right)$ of chain complexes. Since $C\M^k(X)$ is
 an exact category, the agreement theorems apply. Now, by application of the Derived Equivalence Theorem \ref{equivThm}, 
 for a noetherian quasi-projective  scheme $X$, over an affine scheme $\spec{A}$, we obtain zig-zag homotopy equivalences 
 (see (\ref{nota}) for notations)
 \begin{equation}\label{introMainHE}
  \left\{
 \begin{array}{ll}
  {\bf K}\left(C\M^k(X)\right)\iso {\bf K}\left({\SC}h^k\left(\SV(X)\right) \right) &  {\rm of ~the}~{\bf K}{\rm-theory~spaces},\\
 \K\left(C\M^k(X)\right)\iso \K\left({\SC}h^k\left(\SV(X)\right) \right)&  {\rm of ~the}~\K{\rm-theory~spectra}.\\
  \end{array}
 \right.
\end{equation}
Now assume that $X$ is Cohen-Macaulay. In this case, the Thomason-Waldhausen Localization theorem
(\cite[3.2.27]{S1} applies to the inclusion
${\SC}h^{k+1}\left(\SV(X) \right) \hra {\SC}h^k\left(\SV(X)\right)$, and using the identifications (\ref{introMainHE}),
we obtain a  sequence
$$
\diagram 
\K\left(C\M^{k+1}(X)\right) \ar[r] & 
\K\left(C\M^{k}(X)\right) \ar[r] & 
\coprod_{x\in X^{(k)}}\K\left(C\M^{k}(X_x)\right) \\  
\enddiagram 
$$ 
of zig-zag maps (via {\it
homotopy equivalences}) of $\K$-theory spectra, that is a homotopy fibration (see Theorem \ref{Ktheorem}), where $X_x:=\spec{{\CO}_{X,x}}$. 
 This is an
analogue of the homotopy fibration  of $G$-theory spaces, due to Quillen (see proof of \cite[Theorem 5.4]{Q}).
Accordingly, for all integers $n, k\in \Z$ with $k\geq 0$,
there is an exact sequence 
$$
\diagram
\cdots \ar[r] &\K_n\left(C\M^{k+1}(X)\right) \ar[r] & \K_n\left(C\M^{k}(X)\right)\ar[r] &  \oplus_{x\in X^{(k)}}\K_n\left(C\M^{k}(X_x)\right)\\
\ar[r] & \K_{n-1}\left(C\M^{k+1}(X)\right) \ar[r] & \cdots&\\
\enddiagram
$$  
of $\K$-groups (see Corollary \ref{BKiCM}). We remark (\ref{bfkThm}) that, if $X$ is regular,  similar statements regarding ${\bf K}$-theory spaces and groups would also be valid. 
While these results allow us to rewrite the Gersten $\K$-theory complexes in terms of the $\K$-groups of the "local categories"
 $C\M^k(X_x)$ (see Remark \ref{remGerSPSeq}),
they provide further insight in to the same in terms of the $\K$-groups of $C\M^k(X)$. 
When, $A$ is a Cohen-Macaulay local ring with $\dim A=d$ and $X=\spec{A}$, it is a result of Roberts and Srinivas \cite[Proposition 2]{RS}
that the map $K_0\left(C\M^d(X)\right)\iso K_0\left({\SC}h^d\left(\SV(X)\right)\right)$ is an isomorphism, which would be a consequence of the 
above homotopy equivalence (\ref{introMainHE}).

Results on $GW$-theory would be fairly similar. Note that in the diagram (\ref{f5FunctIntroGWadd}) of
equivalences,   the dg categories of $\M^k(X)$ does not have a natural duality structure. Remedy for this was obtained 
by embedding this category in the respective category of perfect complexes.  We assume that $X$ is a quasi-projective scheme over 
an affine scheme $\spec{A}$, with $1/2\in A$. Then, for integers $r=0, 1, 2,3$ and $k\geq 0$, 
we obtain zig-zag homotopy equivalences 
 \begin{equation}\label{introGWHE}
  \left\{
 \begin{array}{ll}
  {\bf GW}^{[r]}\left(dgC\M^k(X)\right)\iso {\bf GW}^{[k+r]}\left(dg^k\SV(X) \right) &  {\rm of ~the}~{\bf GW}{\rm-spectra},\\
{\G}W^{[r]}\left(dgC\M^k(X)\right)\iso {\G}W^{[r]}\left(dg^k\SV(X) \right)&  {\rm of ~the}~{\G}W{\rm-bispectra}.\\
  \end{array}
 \right.
\end{equation}
When $X$ is Cohen-Macaulay,
there is sequence of zig-zag maps
$$
\diagram 
{\G}W^{[-1+r]}\left(dgC\M^{k+1}(X) \right) \ar[r] &
{\G}W^{[r]}\left(dgC\M^{k}(X)  \right) \ar[r] &
\coprod_{x\in X^{(k)}}{\G}W^{[r]}\left(dgC\M^{k}\left(X_x\right)  \right)  \\
\enddiagram
$$
that is a homotopy fibration of ${\G}W$-bispectra. When $X$ is regular, similar homotopy firation of  ${\bf GW}$-spectra would also 
be valid. Again, this is a ${\G}W$-analogue of the homotopy fibration of the $G$-theory spaces of Quillen (see the proof of
\cite[Theorem 5.4]{Q}). Note that such  statements about Grothendieck Witt theory would not make sense in $G$-theory  
(meaning, working with $Coh^k(X)$) because of non-existence of any natural duality in the respective categories. 

With respect to implications to the derived Witt theory, assume  that $X$ is regular quasi-projective 
scheme over an affine scheme $\spec{A}$, with $1/2\in A$. 
For an integer $k\geq 0$,  consider the  exact sequence of the derived categories \cite{BW}:
$$
\diagram 
{\SD}^{k+1}\left(\SV(X)\right) \ar[r] & {\SD}^k\left(\SV(X)\right) \ar[r]& \coprod_{x\in X^{(k)}} {\SD}^k\left(\SV(X)\right). 
\enddiagram
$$  
Then, the twelve term exact sequence of Witt groups, due to Balmer (\cite[Corollary 6.6]{TWG1}), corresponding to this sequence,
 reduces two five term   exact sequences (see Theorem \ref{witt5Term}), one of them being the following:
 $$
\diagram 
0\ar[r] & W^{-1}\left({\CD}^b\left(C\M^{k+1}(X)\right) \right)  \ar[r] & 
W\left(C\M^{k}(X)\right)  \ar[r] &\oplus_{x\in X^{(k)}} W\left(C\M^k(X_x) \right) \\
\ar[r]  & W\left(C\M^{k+1}(X)\right)  \ar[r] &W^{1}\left({\CD}^b\left(C\M^{k}(X)\right) \right) \ar[r]&0\\
\enddiagram 
$$

We  point out that,  
contrary to the usual filtration by co-dimension of the support, in this article, for a scheme $X$ we use the filtrations 
$\M^k(X)\subseteq \M(X)$ and $Coh^k(X) \subseteq Coh(X)$ by grade (see (\ref{nota})). 
When $X$ is Cohen-Macaulay, these filtrations coincide with the   filtration by co-dimension of the support.
This article is  a culmination of an initiative (\cite{MS1, MS2, MNoe, SS, Mf, FM}) 
to place the category $C\M^k(X)$ at a rightful position in the Algebraic $K$-theory, $GW$-theory and Witt theory of schemes
$X$
and the respective Gersten complexes. 
This category $C\M^k(X)$ behaves like the category of modules of finite length and finite projective dimension, at co-dimension $k \leq \dim X$. 
Significant amount of study of the category  $C\M^k(X)$ was done in \cite{MS1, MNoe, Mf}.

Before we close this introduction, we comment on the lay out of this article. 
In \S \ref{foxMORP} we recall or prove some preliminaries that we need. 
The Derived Equivalence Theorem \ref{equivThm} is established in \S \ref{MAIN}.
In \S \ref{SecKTheo} we establish the implications in $K$-theory. 
We deal with $GW$-theory in \S \ref{secGrotWitt}. In \S \ref{secDerWitt} we discuss Derived Witt theory. 
In 
\S  \ref{secSpGWKGW}, we  give some background information on ${\bf GW}$-spectrum and ${\G}W$-bispectrum. 

\section{Preliminaries}\label{foxMORP}

First, we set up some notations. 
\begin{notations}
\label{nota}
{\rm Throughout this article, $X$ will denote a quasi-projective scheme
over a noetherian affine scheme $Spec(A)$, with finite dimesion $d:=\dim X$. We introduce 
further notations. 
\bE
\item For $x\in X$, denote $X_x:=\spec{{\CO}_{X,x}}$. 
\item Throughout,  $Coh(X)$ will denote the category of  coherent $\CO_X$-modules
and ${\SV}(X)$ will denote the category of all locally free 
sheaves on $X$.  
\item ({\it Readers are referred to {\rm (\cite[Def. 3.1]{MNoe})}
 for a definition of resolving subcategories of abelian categories}.) 
For a resolving subcategory $\SA$ of an abelian category $\CC$, 
for objects $\CF\in \CC$, $\PD(\CF)$ will denote the minimum of the length of resolutions of $\CF$ by objects
in $\SA$.

Denote  
$$
{\M}({\SA})=\{{\CF}\in \CC: \PD({\CF})<\infty  \}
$$
With $\CC=Coh(X)$ we denote
$$
 \M(X):=\M(\SV(X))=\{\CF\in Coh(X): \PDV(\CF)<\infty \} 
$$ 
\item In this article, we consider  filtration of
$Coh(X)$ and ${\M}(X)$ 
by grade, as opposed to usual filtration by co-dimension of the support.
\bE
\item Recall, for ${\CF}\in Coh(X)$, 
$grade({\CF}):=\min\{t:{\CE}xt^t({\CF}, {\CO}_X)\neq 0  \}$. We remark that, if $X$ is Cohen-Macaulay, then 
$grade(\CF)=co\dim\left(Supp(\CF)\right)$ (see \cite{Mt}). 
\item For integers $k\geq 0$,
denote 
$$
\left\{
\begin{array}{l}
Coh^k(X):=Coh^k_g(X):=\{{\CF}\in Coh(X):grade({\CF})\geq k \}\\
{{\M}^k(X):=\M}_g^k({\SA}):=    \{{\CF}\in{\M}(X) 
:  grade({\CF})\geq k\}\\
C{\M}^{k}(X):=\{\CF\in \M(X): grade(\CF)=k=\PDV(\CF) \}
\end{array}
\right.
$$ 
So, we have a filtration, by grade
$
{\M}({\SA})={\M}^0({\SA}) \supseteq {\M}^1({\SA}) 
\supseteq \cdots \supseteq {\M}^d({\SA}) \supseteq 0.   
$
We will strictly be using this filtration by grade and the notation without the subscript $g$ will be the norm.
Note that ${\M}^k({\SA})$ is a Serre subcategory

of ${\M}({\SA})$ (meaning, it has the "2 out of 3" property) .  
Clearly, when $X$ is  Cohen-Macaulay, this filtration coincides with
the filtration by co-dimension of the support.
\eE%
\item 
For an exact category $\SE$, $Ch^b({\SE})$  
will denote the category of chain complexes. 
The bounded derived category of 
${\SE}$ will be denoted by ${\CD}^b({\SE})$. 
\item For a complex ${\CF}_{\bullet}\in Ch^b\left(Coh(X)\right)$, the homologies will 
be denoted by ${\CH}_i\left({\CF}_{\bullet}\right)$. 

\item \label{SCHSDk}
Also, for ${\SE}={\SV}(X), {\M}(X)$,
and integers $k \geq 0$,  
$$
\left\{
\begin{array}{l}
{\SC}h^k\left(\SE\right):=\left\{ 
\CF_{\bullet}\in Ch^b\left(\SE\right): \forall ~i~\CH_i\left(\CF_{\bullet}\right) \in Coh^k(X)
\right\}\\
{\SD}^k\left(\SE\right):=\left\{ 
\CF_{\bullet}\in {\CD}^b\left(\SE\right): \forall ~i~\CH_i\left(\CF_{\bullet}\right) \in Coh^k(X)
\right\}\\
\end{array}
\right.
$$ 
 would denote the full subcategory of such objects. 
 ({\it Note the difference between two fonts ${\CD}, {\SD}$}.) We remark:
 \bE
 \item 
 ${\SC}h^k\left(\SE\right)$ is a complicial exact category (see \cite{S1} for definition).
 In fact, $\left({\SC}h^k\left(\SE\right), \SQ\right)$ is a complicial exact category with weak equivalences, where 
 weak equivalences being the set $\SQ$ of all quasi-isomorphisms. 
 \item Also, $\SD^k\left(\SE\right)$ is a derived subcategory of ${\CD}^b\left(\SE\right)$.
\eE

\eE  
}  
\end{notations}

We recall the following lemma from \cite{FM}. 
\bL
\label{findLIC} 
Suppose $X$ is a quasi-projective noetherian scheme over $\spec{A}$, with $\dim X=d$. 
Then, 
$X$ is an open subset of $\tilde{X}:=Proj(S)$, for some noetherian
graded ring $S=\oplus_{i=0}^{\infty} S_i$, with $S_0=A$.

Let $Y\subseteq X$ be a closed subset of $X$, with  
 $grade\left({\CO}_Y\right) \geq k$. 
Let $V(I)=\overline{Y}$ be the closure of $Y$,
where $I$ is 
the homogeneous ideal of $S$, defining $\overline{Y}$.
Then,
there is a sequence of  homogeneous elements 
$f_1, \ldots, f_k\in I$ such that $f_{i_1}, \ldots, f_{i_j}$ 
induce regular $S_{({\wp})}$-sequences
$\forall~\wp\in Y\subseteq  X$, and $\forall~
1\leq i_1<i_2< \cdots <i_j\leq k$.  
In particular, \\with $Z=V(f_1, \ldots, f_k) \cap X$, and $\CF_n=\CO_Z^n$, we have 
\bE
\item $\CF_n\in C\M^k(X)$. In fact, $\oplus_{i=1}^n \CO_Z \otimes \CL_i\in C\M^k(X)$ for any locally free sheaves $\CL_i$ of rank one.
\item Further, if $\CG\in Coh^k(X)$, with $Y=Supp\left(\CG\right)$ and  $Z$ as above,  there is a surjective map 
$\CF\sur \CG$ where $\CF:=\oplus_{i=1}^n \CO_Z \otimes \CO_x(n_i)$, for some integers $n_i$.
Note that $\CF\in C\M^k(X)$. 
\eE 
\eL

\begin{remark}\label{choice}
{\rm
In this, one had  
 the choices of the  sequence $f_1, \ldots, f_k$, as required above. 
 We may exploit this flexibility later. 
}\end{remark}

%
%


\bL
\label{projDrop}
Suppose $X$ is a quasi-projective noetherian scheme over $\spec{A}$, with $\dim X=d$. 
 Consider an exact sequence
$$
\diagram 
0\ar[r] & \CK \ar[r] & \CE \ar[r] & \CF \ar[r] & 0\\
\enddiagram
~~{\rm where}~~ \CF\in \M^{k}(X), ~\CE \in C\M^{k}(X).
$$
Then $\PDV(\CK) \leq \max\{k, \PDV(\CF)-1\}$. 
In fact,\\
 $\PDV(\CF)\geq k+1 \Lra \PDV(\CK)=\PD(\CF)-1$. 
\eL
\pf If $\PDV(\CF)=k$ then there is nothing to prove. Assume, $\PDV(\CF)=m\geq k+1$. 
Arguing locally, a simple Tor-argument establishes the lemma.
$\eop$
\bL
\label{finiCMres}
Suppose $X$ is a quasi-projective noetherian scheme over $\spec{A}$, with $\dim X=d$. 
Suppose $\CF\in \M^{k}(X))$.
Then, there is a resolution 
\begin{equation}\label{aCMkresDia}
\diagram 
0\ar[r] & \CE_n\ar[r] & \CE_{n-1} \ar[r]^{\partial_{n-1}} & \cdots \ar[r] & \CE_1 \ar[r] & \CE_0 \ar[r] & \CF \ar[r] & 0\\
\enddiagram 
~~~{\rm with }~\CE_i\in C\M^{k}(X). 
\end{equation}
In fact, $n=\PDV(\CF)-k$. 
\eL 
\pf 
By Lemma \ref{findLIC},  there is a surjective map $\partial_0:\CE_0\sur \CF$, where $\CE_0\in C\M^k(X)$. 
Now, let $\CF_0=\ker(\partial_0)$. If $\PDV(\CF)\geq k+1$, then by lemma \ref{projDrop}, $\PDV(\CF_0)=\PDV(\CF)-1$.  
By repeating this process, we get an exact sequence, as in diagram \ref{aCMkresDia}, with $\CE_n=\ker(\partial_{n-1})$
and $\PDV(\CE_n)=k$. Since $grade(\CE_n)\geq k$, it follows $\CE_n \in C\M^{k}(X)$. 
 \pic $\eop$

\bP
\label{resPro}{\rm 
Suppose $X$ is a quasi-projective noetherian scheme over $\spec{A}$, with $\dim X=d$. 
 Suppose 
$$
\diagram 
0\ar[r] & \CK \ar[r] & \CF \ar[r] & \CG\ar[r] & 0 \\
\enddiagram 
\qquad
{\rm be~exact~in}~\qquad \M^{k}(X)
$$
Then, 
\bE
\item First,
$$
\CK, \CG \in C\M^{k}(X) \Lra ~~~ \CF \in C\M^{k}(X) 
$$
\item Then,
$$
\CF, \CG \in C\M^{k}(X) \Lra ~~~ \CK \in C\M^{k}(X) 
$$
\eE
}\eP
\pf The proof follows by routine chasing the long exact sequence of the ${\CE}xt$-modules.

$\eop$

\bC
\label{itIsRes}
Let $X$ be a quasi-projective scheme over an affine scheme $\spec{A}$. Then, $C\M^{k}(X) \subseteq Coh^{k}(X)$ is a resolving subcategory.
Further, 
$$
\M^{k}(X)=\{\CF\in Coh^k(X): \dim_{C\M^{k}(X)}(\CF)<\infty \}
$$
\eC
\pf 
Let $\CG\in Coh^k(X)$ be an object. By Lemma \ref{findLIC}, there is a surjective morphism $\CF\sur \CG$ with $\CF\in C\M^k(X)$.
Now it follows from Lemma \ref{resPro} that $C\M^k(X)$ is a resolving subcategory of $Coh^k(X)$.

 Now, suppose $\CF\in \M^k(X)$. By Lemma \ref{finiCMres}, $\dim_{C\M^{k}(X)}(\CF)<\infty$. Conversely, suppose $\CF\in Coh^k(X)$
 and
$\dim_{C\M^{k}(X)}(\CF)<\infty$. In particular, $\PDV(\CF)<\infty$. So, $\CF\in \M(X)$ and hence $\CF\in \M^k(X)$.
  \pic $\eop$ 
 
\subsection{Some further preparation} 
In this subsection, we gather some preparatory information that will be used more directly in the proof of the main theorems. 
Among them, the following lemma would be of some interest to us.
\bP
\label{dominate} 
Suppose $\SA$ is a resolving subcategory of an abelian category $\SC$. Suppose 
$$
\diagram 
0\ar[r] & \CG_n\ar[r]_{\partial_n}  & \CG_{n-1} \ar[r]_{\partial_{n-1}}  & \cdots \ar[r] & \CG_1 \ar[r]_{\partial_1}  & \CG_0 \ar[r]_{\partial_0} &\CG\ar[r] & 0\\ 
\enddiagram
$$
is an exact sequence in $\SC$ and $f:\CF\to \CG$ is a morphism in $\SC$. Assume 
$\dim_{\SA}(\CF)<\infty$ and every $\SA$-resolution of $\CF$ terminates. Then, there is 
a complex $\CE_{\bullet}\in Ch^b(\SA)$ such that
\bE
\item $\CE_i=0~\forall~i\leq -1$.
\item $\CE_{\bullet}$ is a resolution of $\CF$.
\item There is a morphism $f_{\bullet}:\CE_{\bullet}\lra \CG_{\bullet}$ of complexes, such that $\CH_0(f_{\bullet})=f$.
\eE
\eP 
\pf Consider the diagram
$$
\diagram 
\CE_0\ar@{->>}[r]^{\partial} &\Gamma \ar[r] \ar[d] & \CF\ar[d]^f\\
&\CG_0\ar@{->>}[r] & \CG\\ 
\enddiagram 
$$
where $\Gamma$ is the pullback, $\CE_0\in \SA$ and $\CE_0\sur \Gamma$ is a surjective map. 
Note that all horizontal arrows are surjective. Let $d_0:\CE_0\sur \CF$ be the composition. 
Now, we use induction. Suppose, $\CE_m$ and $f_m: \CE_m \to \CG_m$ were defined. Consider the diagram: 
$$
\diagram 
\CE_{m+1}\ar@{->>}[r]^{\epsilon} &\Gamma_m\ar@{->>}[r]^p\ar[d]&Z\ar[r]\ar[d]&\CE_m\ar[r]^{d_m}\ar[d]^{f_m} & \CE_{m-1}\ar[d]\\
&\CG_{m+1}\ar@{->>}[r]_{\partial_{m+1}}&B\ar[r]&\CG_m \ar[r]_{\partial_r} & \CG_{m-1} \\ 
\enddiagram 
$$
where $B=\ker(\partial_m)$ and $Z=\ker(d_r)$, $\Gamma_m$ is the pullback, $\CE_{m+1}\in \SA$ and $\epsilon$ is a surjective map. Now, let 
$d_{m+1}:\CE_{m+1}\lra \CE_m$ and $f_{m+1}: \CE_{m+1}\to \CG_{m+1}$
be the respective compositions. Since every resolution of $\CF$ terminates, this process terminates, at least when $m\gg n$. 
\pic $\eop$.

\bC
\label{forCAisCMk} 
Suppose $X$ is a quasi-projective noetherian scheme over an affine scheme $\spec{A}$, with $\dim X=d$. 
Consider the diagram
$$
\diagram 
&&&&&&\CF\ar[d]^f&\\
0\ar[r] & \CG_n\ar[r] & \CG_{n-1} \ar[r] & \cdots \ar[r] & \CG_1 \ar[r] & \CG_0 \ar[r] &\CG\ar[r] & 0\\ 
\enddiagram
$$
where $\CF\in \M^k(X)$ and the second line is an exact sequence in $Coh^k(X)$. Then, 
\bE
\item There is complex $\CE_{\bullet}\in Ch^b(C\M^k(X))$, with $\CE_i=0~\forall~i\leq -1$,
\item $\CE_{\bullet}$ is a resolution of $\CF$,
\item there is a morphism $f_{\bullet}: \CE_{\bullet} \lra \CG_{\bullet}$
of complexes, such that $\CH_0(f_{\bullet})=f$. 
\eE
\eC
\pf Note by Lemma \ref{projDrop}, every $C\M^k(X)$-resolution of $\CF$ terminates. Now, the proof is complete by 
Proposition \ref{dominate}. $\eop$

The following   version of ( \cite[Lemma 3.3]{FM})
would  be more precise for our purpose. 
\bL\label{1TONHom}
Let ${\SC}$ be  an abelian category and $\SA$ be an exact subcategory, such that
all epimorphisms in $\SA$ are admissible.  Let ${\CF}_{\bullet},
{\CG}_{\bullet}$ be two objects in ${\CD}^b({\SA})$.  
Assume (1) ${\CH}_r({\CF}_{\bullet})=0~\forall r\leq -1$,
(2) ${\CH}_r({\CG}_{\bullet})=0~\forall r\geq 0$.
Then, $Mor_{{\CD}^b({\SA})}({\CF}_{\bullet},
{\CG}_{\bullet})=0$, if $\ker\left(\CG_0\to \CG_{-1}\right) \in \SA$.
This would be the case in the following situations (1) $\CG_{\bullet}$ is concentrated at a single degree $n_0\leq -1$, 
(2) $\forall~i\leq -1,~\CH_i\left(\CG_{\bullet}\right)\in \SA$, or (3) $\SA$ has the 2 out of 3 property. 
\eL
\pf 
Let $\eta_{\bullet}:{\CF}_{\bullet}
\lra {\CG}_{\bullet}$ be a morphism.
We can assume, $\eta_{\bullet}$ is a map of complexes (denominator
free). Further, by replacing by a quasi-isomorphic complex, 
we assume ${\CF}_i=0~\forall~i\leq -1$. 
Define the subcomplex ${\CG}_{\bullet}'\hookrightarrow {\CG}_{\bullet}$ 
by setting ${\CG}_i'={\CG}_i ~\forall ~i\geq 1$, ${\CG}_0'=\ker({\CG}_0
\lra {\CG}_{-1})$ and ${\CG}_i'=0 ~\forall ~i\leq -1$. 
Note that  ${\CG}'_{\bullet}$ is a complex in $Ch^b(\SA)$ and is exact. Hence, ${\CG}_{\bullet}'\cong 0$ in ${\CD}^b({\SA})$. 
Since, $\eta_{\bullet}$ factors through ${\CG}_{\bullet}'$, $\eta_{\bullet}=0$. 
\pic$\eop$

The following 
 special case of the faithfulness property is a variation of the same in \cite{FM, SS}. 
\bL\label{1tonFaith}
Let ${\SC}$ be  an abelian category and $\SA\hra \SB$ is two exact subcategories of $\SC$, such that
all epimorphisms in $\SA$  are admissible. Let $\Phi: {\CD}^b\left(\SA\right) \lra {\CD}^b\left(\SB\right)$
denote the induced functor of the bounded derived categories. 
Suppose $f_{\bullet}:{\CF}_{\bullet}\lra {\CG}_{\bullet}$ is a morphism in 
${\CD}^b(\SA)$ such that 
(a) ${\CH}_i({\CF}_{\bullet})=0~\forall~i\leq n_0-1$,
(b) 
$\CG_i=0~\forall~i\neq n_0$, and 
(c)  $\Phi\left(f_{\bullet}\right)=0$ in ${\CD}^b(\SB)$. 
Then, $f_{\bullet}=0$ in ${\CD}^b(\SA)$. 
\eL
\pf 
First, write $f_{\bullet}=g_{\bullet}
t_{\bullet}^{-1}$ where $t_{\bullet}$ is a quasi-isomorphism
and $g_{\bullet}$ is a chain complex map. By usual arguments, we can assume that $f_{\bullet}=g_{\bullet}$ is a chain complex map
and $n_0=0$.
So, $\forall~i\leq -1,~\CH_i(\CF_{\bullet})=0$ and 
$\forall~i\neq 0~\CG_i=0$. 
Since, all epimorphisms in $\SA$ are admissible, we can assume $\forall~i\leq -1,~\CF_i=0$. 

Since, $\Phi(f_{\bullet})=0$, it follows $\CH_0\left(f_{\bullet}\right)=0: \CH_0\left(\CF_{\bullet}\right) \lra \CG_0$. 
With $p:{\CF}_0\lra {\CH}_0({\CF}_{\bullet})$, we have $f_0
={\CH}_0(f_{\bullet}) p=0$.  
Hence $f_{\bullet}=0$. 
\pic $\eop$

For the convenience of the readers we include the following well known lemma.
\bL\label{zeroIn3angle}
Suppose $\ST$ is a triangulated category and 
$
\diagram 
T^{-1}\Delta \ar[r]^t & \CF \ar[r]^f & \CG \ar[r]^h & \Delta
\enddiagram
$
is an exact triangle. Then, $f=0$ if and only if $t$ splits. 
\eL

\section{The Equivalence Theorems}\label{MAIN}
In this section we state and prove the  main equivalence theorems. 
\bT
\label{mainThmFoxII}
Let $X$ be a noetherian quasi-projective scheme over an affine scheme $\spec{A}$
and $k\geq 0$ be a fixed  integer. Consider
the inclusion functor 
$C\M^{k}(X) \hra \M^{k}(X)$ and let
$$
 \zeta: \CD^b\left(C\M^{k}(X)\right) \lra
\CD^b\left(\M^{k}(X)\right) 
\quad {\rm denote~the ~induced~functor}
$$
of the derived categories. Then $\zeta$ is an equivalence of derived categories.
\eT
We will split this proof in to three propositions as follows. Before that, we introduce the following definition for subsequent reference.
\bD
Suppose $X$ is a noetherian scheme. 
For a complex ${\CF}_{\bullet}\in Ch^b(Coh(X))$,
we say that $width({\CF}_{\bullet})\leq r$, if ${\CH}_i({\CF}_{\bullet})=0$ unless $n\geq i \geq n-r$, for some integer $n$. 
\eD

\bP
\label{Mkp1toMkfull}
Let $X$ be a noetherian quasi-projective scheme as in (\ref{mainThmFoxII}) and $k\geq 0$ be a fixed  integer. 
Then
the functor $\zeta$ is  full. 
\eP
\pf 
Suppose $\CF_{\bullet}, \CG_{\bullet}$ are two complexes in $Ch^b\left(C\M^k(X) \right)$. We will prove that the map
$$
\zeta: Mor_{{\CD}^b(C\M^{k}(X))}\left( {\CF}_{\bullet}, 
{\CG}_{\bullet}\right) 
\to Mor_{{\CD}^b(\M^{k}(X))}\left( {\CF}_{\bullet}, 
{\CG}_{\bullet}\right) 
\quad {\rm is~surjective}. 
$$
 Without loss of generality, we assume that $\forall ~i\leq -1~\CH_i\left(\CF_{\bullet}\right)=\CH_i\left(\CG_{\bullet}\right)=0$.
Since, $C\M^k(X)$ is a resolving subcategory, 
using the usual arguments we can assume that $\forall~i\leq -1~\CF_i=\CG_i=0$.
Now, suppose $f_{\bullet}:{\CF}_{\bullet}
\to {\CG}_{\bullet}$ is a morphism in ${\CD}^b({\M}^{k}(X))$.  We will prove that $f_{\bullet}$ is in the image of $\zeta$.
We will prove this by induction 
on $r:= width({\CF}_{\bullet}\oplus {\CG}_{\bullet})$.

Let $r=0$.
We have,
$f_{\bullet}=g_{\bullet}t_{\bullet}^{-1}:
\diagram {\CF}_{\bullet} & W_{\bullet}\ar[r]^{g_{\bullet}}\ar[l]_{t_{\bullet}}
& {\CG}_{\bullet}  \enddiagram$, where $W_{\bullet}$ is in $\CD^b(\M^k(X))$ and $t_{\bullet}$ is a quasi-isomorphism.
By same argument, we assume $\forall~i\leq -1,~W_i=0$. 
By Corollary \ref{forCAisCMk}, there is a quasi-isomorphism $\epsilon_{\bullet}: \CE_{\bullet}\lra W_{\bullet}$, where 
$\CE_{\bullet}\in Ch^b\left(C\M^k(X)\right)$.
Therefore, $f=g_{\bullet}t_{\bullet}^{-1}= \zeta\left(\left(g_{\bullet}\epsilon_{\bullet}\right)\left(t_{\bullet}\epsilon_{\bullet}\right)^{-1}\right)$. 

Now suppose $r>0$ and  $width({\CF}_{\bullet}\oplus 
{\CG}_{\bullet})=r$.
As before, we have 
${\CF}_i={\CG}_i=0$
for all $i\leq -1$. Further, assume that  
either ${\CH}_0({\CF}_{\bullet})\neq 0$ or ${\CH}_0({\CG}_{\bullet})\neq 0$.
As before,   
$f_{\bullet}:{\CF}_{\bullet}\lra {\CG}_{\bullet}$ is a morphism 
in ${\CD}^b({\M}^k(X))$.
By (\ref{itIsRes}), 
there are objects
${\CE}, {\CL} \in C\M^{k}(X)$ 
and  surjective morphisms $\CE \sur \CF_0$, $\CL\sur \CG_0$.  
Let $\CE_{\bullet}$ denote the complex with $\CE_0=\CE$ and $\CE_i=0~\forall ~i\neq 0$ and likewise,
let $\CL_{\bullet}$ denote the complex defined by $\CL$, concentrated at degree zero. 
Let $\nu:{\CE}_{\bullet}\lra {\CF}_{\bullet}$, 
$\mu:{\CL}_{\bullet}\lra {\CG}_{\bullet}$ denote the obvious morphisms.
Then, 
\bE
\item Note $\CH_0(\CE_{\bullet})=\CE$ and $\CH_0(\CL_{\bullet})=\CL$. 

\item  
By replacing ${\CL}_{\bullet}$ by ${\CE}_{\bullet}\oplus {\CL}_{\bullet}$,
we can further assume that the  diagram
\begin{equation}\label{inAbove}
\diagram
{\CE}_{\bullet} \ar[rr]^{\nu} \ar[d]_{\varphi} && 
{\CF}_{\bullet}\ar[d]^f\\
{\CL}_{\bullet} \ar[rr]_{\mu} & & {\CG}_{\bullet}\\
\enddiagram 
\qquad {\rm commutes}.
\end{equation} 
\eE 
Embed $\mu, \nu$  in  exact triangles  and 
obtain the vertical morphism of exact triangles:
\begin{equation}\label{Jul6gotDia9J}
\diagram
T^{-1}\Delta_{\bullet}\ar[r]^{\nu_0}\ar[d]_{T^{-1}\eta} 
& {\CE}_{\bullet}  \ar[r]^{\nu} \ar[d]_{\varphi} & 
{\CF}_{\bullet}\ar[d]^f\ar[r]^{\nu_1} & \Delta_{\bullet}\ar[d]^\eta\\
T^{-1}\Gamma_{\bullet} \ar[r]_{\mu_0}
&{\CL}_{\bullet} \ar[r]_{\mu}  & {\CG}_{\bullet} \ar[r]_{\mu_1}
& \Gamma_{\bullet}\\ 
\enddiagram
\end{equation}
where the two horizontal exact triangles are in $Ch^b(C{\M}^{k}(X))$   
and the vertical morphisms are in ${\CD}^b({\M}^k(X))$.
Since $width(\Delta_{\bullet}\oplus \Gamma_{\bullet})\leq r-1$ and $width(\CE_{\bullet}\oplus \CL_{\bullet})=0$, 
the induction hypotheses applies to $\eta, \varphi$. So, there 
are morphisms 
$\tilde{\eta}: \Delta_{\bullet}\lra \Gamma_{\bullet}$,
$\tilde{\varphi}: {\CE}_{\bullet} \lra {\CL}_{\bullet}$
in ${\CD}^b(C\M^k(X))$ 
such that $\zeta(\tilde{\eta})=\eta$, $\zeta(\tilde{\varphi})=\varphi$ ({\it in fact, by construction}, $\tilde{\varphi}=\varphi$). 
This gives the following commutative (as clarified below) diagram:
\begin{equation}\label{1TONBelow9J}
\diagram
T^{-1}\Delta_{\bullet}\ar[r]^{\nu_0}\ar[d]_{T^{-1}\tilde{\eta}} 
& {\CE}_{\bullet} \ar[r]^{\nu} \ar[d]_{\tilde{\varphi}} & 
{\CF}_{\bullet}\ar@{-->}[d]^{g}\ar[r]^{\nu_1} 
& \Delta_{\bullet}\ar[d]^{\tilde{\eta}}\\
T^{-1}\Gamma_{\bullet} \ar[r]_{\mu_0}
&{\CL}_{\bullet} \ar[r]_{\mu}  & {\CG}_{\bullet} \ar[r]_{\mu_1}
& \Gamma_{\bullet}\\ 
\enddiagram
\qquad {\rm in}~~{\CD}^b(C{\M}^{k}(X)). 
\end{equation}
Since the image of the left square in ${\CD}^b({\M}^k(X))$
commutes, it follows from  Lemma \ref{1tonFaith},  that the left square commutes in ${\CD}^b\left(C\M^{k}(X)\right)$.
Therefore, by properties of triangulated categories, there  is a morphism $g$ in ${\CD}^b(C{\M}^{k}(X))$,
making the diagram commute in ${\CD}^b\left(C\M^k(X)\right)$. Apply  
$\zeta$ to (\ref{1TONBelow9J}) and compare with (\ref{Jul6gotDia9J}). 
With $h=f-\zeta(g)$ we obtain the commutative diagram
\\
$$
\diagram
T^{-1}\Delta_{\bullet}\ar[r]^{\nu_0}\ar[d]_{0} 
& {\CE}_{\bullet} \ar[r]^{\nu} \ar[d]_{0} & 
{\CF}_{\bullet}\ar[d]_{h}\ar@{-->}[dl]_{\beta'}\ar[r]^{\nu_1} 
& \Delta_{\bullet}\ar[d]^{0}\ar@{-->}[dl]^{\epsilon}\\
T^{-1}\Gamma_{\bullet} \ar[r]_{\mu_0}
&{\CL}_{\bullet} \ar[r]_{\mu}  & {\CG}_{\bullet} \ar[r]_{\mu_1}
& \Gamma_{\bullet}\\ 
\enddiagram
\qquad {\rm in}~~{\CD}^b({\M}^k(X)), 
$$
where $\beta'$ and $\epsilon$ are given by weak kernel and weak cokernel
proporties.
Consider the case ${\CH}_0({\CF}_{\bullet})=0$. Then, by lemma \ref{1TONHom}, 
$Mor({\CF}_{\bullet}, {\CL}_{\bullet})=0$.  
Therefore, 
$\beta'=0$. Hence $h=0$ and $f=\zeta(g)$. 
So, it is established, whenever 
${\CH}_0({\CF}_{\bullet})=0$ and $width({\CF}_{\bullet} \oplus 
{\CG}_{\bullet})\leq r$, then 
$$
Mor_{{\CD}^b(C{\M}^{k}(X))}
({\CF}_{\bullet}, {\CG}_{\bullet})\sur 
Mor_{{\CD}^b({\M}^{k}(X))}
({\CF}_{\bullet}, {\CG}_{\bullet})
\quad {\rm is~ surjective.}
$$  
This fact will be exploited  in the next step. 
Now, assume ${\CH}_0({\CF_{\bullet}})\neq 0$.
It follows, $width(\Delta_{\bullet}\oplus  
{\CG}_{\bullet})\leq r$ and ${\CH}_0({\Delta_{\bullet}})=0$. 
Therefore, the map\\
$Mor_{{\CD}^b(C{\M}^{k}(X))}
({\Delta}_{\bullet}, {\CG}_{\bullet})\sur 
Mor_{{\CD}^b({\M}^k(X))}
({\Delta}_{\bullet}, {\CG}_{\bullet})$ is surjective.
Hence, $\zeta(\tilde{\epsilon})=\epsilon$ for some 
$\tilde{\epsilon}\in 
Mor_{{\CD}^b(C{\M}^{k}(X))}
({\Delta}_{\bullet}, {\CG}_{\bullet})$.
So, $h=\zeta(\tilde{\epsilon}) \nu_1=\zeta(\tilde{\epsilon} \nu_1)$ is in the image of $Mor_{{\CD}^b(C{\M}^{k}(X))}
({\CF}_{\bullet}, {\CG}_{\bullet})$. Hence, so is $f=\zeta(g)+h$. 
\pic
$\eop$ 

The following establishes faithfulness of $\zeta$. 
\bP
\label{Mkp12Mkfaith}
Let $X$ be a noetherian quasi-projective scheme as in (\ref{mainThmFoxII}) and $k\geq 0$ be a fixed integer. 
Then,
the functor $\zeta: {\CD}^b\left(C\M^{k}(X)\right)\lra {\CD}^b\left(\M^k(X)\right)$ is faithful.
\eP 
\pf Let $f_{\bullet}:{\CF}_{\bullet}\lra {\CG}_{\bullet}$
be a morphism in ${\CD}^b\left(C{\M}^{k}(X)\right)$ such that
$\zeta(f_{\bullet})=0$. We need to prove that  $f_{\bullet}=0$ in ${\CD}^b\left(C{\M}^{k}(X)\right)$. Without
loss of any generality, we can assume that 
$f_{\bullet}$ is in $Ch^b({C\M}^k(X))$. 
Embed $f_{\bullet}$ is an exact triangle
$$
\diagram
T^{-1}\Delta_{\bullet} \ar[r]^g & {\CF}_{\bullet} \ar[r]^{f_{\bullet}} & {\CG}_{\bullet} 
\ar[r]^h & \Delta_{\bullet}\\ 
\enddiagram
\quad 
{\rm in}~~{\CD}^b(C{\M}^{k}(X))
$$ 
Under $\zeta$, this triangle maps to 
$
\diagram
T^{-1}\Delta_{\bullet} \ar[r]^g & {\CF}_{\bullet} \ar[r]^0 & {\CG}_{\bullet} 
\ar[r]^h & \Delta_{\bullet}\\ 
\enddiagram
$
in 
${\CD}^b({\M}^k(X)).
$
Therefore, by Lemma
\ref{zeroIn3angle}, 
there is a spilt $\eta:{\CF}_{\bullet}\lra T^{-1}\Delta_{\bullet}$,
in ${\CD}^b({\M}^k(X))$, 
of $g$. So, $g\eta=1_{{\CF}_{\bullet}}$ in ${\CD}^b({\M}^k(X))$. 
Since $\zeta$ is full  (\ref{Mkp1toMkfull}), we can assume that $\eta$ is in 
${\CD}^b(C{\M}^{k}(X))$. Now, embed $g\eta$ in a triangle:
$$
\diagram
{\CF}_{\bullet}\ar[r]^{g\eta} & {\CF}_{\bullet} \ar[r]& 
\Gamma_{\bullet}\ar[r] & T{\CF}_{\bullet} 
\enddiagram
\quad {\rm in}~~{\CD}^b(C{\M}^{k}(X))
$$
This triangle maps to 
$
\diagram
{\CF}_{\bullet}\ar@{=}[r] & {\CF}_{\bullet} \ar[r]& 
\Gamma_{\bullet}\ar[r] & T{\CF}_{\bullet} 
\enddiagram
$
in 
${\CD}^b({\M}^k(X))$. 
Therefore, $\Gamma_{\bullet}=0$ in ${\CD}^b({\M}^k(X))$.
That means, $\Gamma_{\bullet}$ is exact, as a complex.
Hence, $\Gamma_{\bullet}\cong 0$ in ${\CD}^b(C{\M}^{k}(X))$. 
So, $g\eta$ is an isomorphism
in ${\CD}^b(C{\M}^k(X))$.  
Therefore $g(\eta(g\eta)^{-1})=1$ and hence 
$\eta(g\eta)^{-1}$ is a split of $g$ 
in ${\CD}^b(C{\M}^{k}(X))$. 
Therefore, by Lemma \ref{zeroIn3angle}, 
$f_{\bullet}=0$ in ${\CD}^b(C{\M}^{k}(X))$.
\pic $\eop$ 

The third of the three part proof 
of  Theorem  \ref{mainThmFoxII}, would be the following proposition on essential surjectivity of $\zeta$.

\bP
\label{essSurjective}
Let $X$ be a noetherian quasi-projective scheme as in (\ref{mainThmFoxII}) and $k\geq 0$ be a fixed  integer. 
Then,
the functor $\zeta: {\CD}^b\left(C\M^{k}(X)\right)\lra {\CD}^b\left(\M^k(X)\right)$ is essentially surjective.
\eP
\pf Let $\CF_{\bullet}\in {\CD}^b\left(\M^k(X)\right)$ be an object. We will prove that there is a complex $\tilde{\CF}_{\bullet}\in 
{\CD}^b\left(C\M^k(X)\right)$ such that $\CF_{\bullet}\ \cong \zeta\left(\tilde{\CF}_{\bullet}\right)$. We will use induction 
on $r:=width\left({\CF}_{\bullet}\right)$.
As usual, we can assume that $\forall~i\leq -1~\CF_i=0$
and $\CH_0\left(\CF_{\bullet}\right)\neq 0$. 

Suppose $r=0$.  In this case, $\forall~i\neq 0~\CH_i\left(\CF_{\bullet}\right)=0$.
It follows that $\CH_0\left(\CF_{\bullet}\right) \in \M^k(X)$. In fact, $\CF_{\bullet}$ is quasi-isomorphic to the complex,
concentrated at degree zero, defined by $\CH_0\left(\CF_{\bullet}\right)$. Now, it follows from Lemma \ref{finiCMres}, 
that such a 
complex is quasi-isomorphic to a complex $\CE_{\bullet} \in C\M^k(X)$. This settles the base case, $r=0$. 

Now, assume $r\geq 1$. 
By (\ref{itIsRes}), there is an epimorphism $\CE\sur \CF_0$, with $\CE\in C\M^k(X)$. 
 Let ${\CE}_{\bullet}$ denote the complex with $\CE_0=\CE$ and $\CE_i=0~\forall ~i\neq 0$
 and let $\nu: \CE_{\bullet} \lra \CF_{\bullet}$ denote the obvious map. Then, 
 $\CH_0(\nu): \CH_0(\CE_{\bullet}) \sur \CH_0(\CF_{\bullet})$ is surjective. 
 
Embed $\nu$ in an exact triangle in $\CD^b({\M}^k(X))$: 
$
\diagram
T^{-1}{\Delta}_{\bullet} \ar[r]^{\nu_0}& 
{\CE}_{\bullet}\ar[r]^{\nu} & {\CF}_{\bullet}\ar[r]^{\nu_2} 
& {\Delta}_{\bullet}.  
\enddiagram 
$ 
It follows from the long exact sequence of the homologies
$$
{\CH}_i({\Delta}_{\bullet})=\left\{
\begin{array}{ll}
0 & if ~i\leq 0\\
{\CH}_i({\CF}_{\bullet}) & if ~i\geq 2\\
\end{array}
\right.
$$
and $
\diagram
0\ar[r] & {\CH}_1({\CF}_{\bullet}) 
\ar[r] & {\CH}_1({\Delta}_{\bullet}) \ar[r] &
{\CH}_0({\CE}_{\bullet})\ar[r] & {\CH}_0({\CF}_{\bullet})
\ar[r] & 0\\ 
\enddiagram
$
is exact.\\
Therefore, $width({\Delta}_{\bullet})\leq r-1$. By induction, there is a 
complex $\tilde{\Delta}_{\bullet}\in {\CD}^b(C{\M}^k(X))$ such that
$\zeta(\tilde{\Delta}_{\bullet})\cong {\Delta}_{\bullet}$.  
Since $\zeta$ is full (\ref{Mkp1toMkfull}), there 
is a morphism $\eta_0:T^{-1}\tilde{\Delta}_{\bullet}\lra 
{\CE}_{\bullet}$ in ${\CD}^b(C{\M}^k(X))$ such that
$\zeta(\eta_0)=\nu_0$. Now, embed $\eta_0$ in
an exact triangle in ${\CD}^b(C{\M}^k(X))$:
$$
\diagram
T^{-1}\tilde{{\Delta}}_{\bullet} \ar[r]^{\eta_0}& 
{\CE}_{\bullet} \ar[r]^{\eta} & {\CU}_{\bullet}\ar[r]^{\eta_2} & 
\tilde{{\Delta}}_{\bullet}\\  
\enddiagram. 
$$
Now apply  $\zeta$ to this triangle and complete the diagram:
$$
\diagram
T^{-1}{\Delta}_{\bullet}\ar[d]_{\wr}\ar[r]^{\nu_0}& 
{\CE}_{\bullet}\ar[r]^{\nu}\ar@{=}[d] 
& {\CF}_{\bullet}\ar@{-->}[d]_{\wr}^{\epsilon}\ar[r]^{\nu_2} 
& {\Delta}_{\bullet} \ar[d]^{\wr}\\
T^{-1}\zeta\left(\tilde{\Delta}_{\bullet}\right) \ar[r]_{~~~\zeta(\eta_0)}& 
{\CE}_{\bullet} \ar[r]_{\zeta(\eta)} 
& \zeta({\CU}_{\bullet}) \ar[r]_{\zeta(\eta_2)} & 
\zeta\left(\tilde{\Delta}_{\bullet}\right)\\  
\enddiagram 
\qquad {\rm in} \quad {\CD}^b(\M^k(X). 
$$
The isomorphism $\epsilon \in {\CD}^b(\M^k(X)$ is obtained by properties of 
triangulated categories. This completes the proof of the proposition. $\eop$

\vspace{5mm}
Finally, we are ready to formally complete the proof of Theorem \ref{mainThmFoxII}. 

\vspace{5mm} 
\noindent{\bf The  proof of Theorem \ref{mainThmFoxII}}: The proof is complete by Propositions 
 \ref{Mkp1toMkfull}, \ref{Mkp12Mkfaith} and \ref{essSurjective}.  $\eop$
 

Using the same proof as above we obtain the following. 
\bT
\label{mkp1toMkFlF}
Let $X$ be a noetherian quasi-projective scheme as in (\ref{mainThmFoxII}) and $k\geq 0$ be a fixed  integer. 
Consider the inclusion functor $\M^{k+1}(X) \to \M^k(X)$.
Then, the induced functor $\beta: {\CD}^b\left(\M^{k+1}(X) \right) \lra  {\CD}^b\left(\M^{k}(X)\right)$ is fully faithful. 
Consequently, so is the functor ${\CD}^b\left(C\M^{k+1}(X)\right)\lra  {\CD}^b\left(\M^{k}(X)\right)$. 
\eT 
\pf The proof of fullness and faithfulness of $\beta$ are, respectively, are similar to that of
Propositions  \ref{Mkp1toMkfull} and \ref{Mkp12Mkfaith}. The second statement 
follows from the former, because the functor $\zeta: {\CD}^b\left(C\M^{k+1}(X) \right)
\iso {\CD}^b\left(\M^{k+1}(X) \right)$ is an equivalence, by  Theorem \ref{mainThmFoxII}.
\pic $\eop$

Combining the results in \cite{FM}, we summarize the results as follows. 
\bT
\label{equivThm}
Let $X$ be a noetherian quasi-projective scheme as in (\ref{mainThmFoxII}) and $k\geq 0$ be a fixed  integer. 
Consider the commutative diagram of  functors of derived categories:
\begin{equation}\label{fiveFunct}
\diagram
{\CD}^b(C\M^{k+1}(X))\ar[r]^{\zeta}_{\sim}\ar[d]^{\alpha} &{\CD}^b({\M}^{k+1}(X))\ar[r]^{\iota}_{\sim}\ar[d]^{\beta}
& {\SD}^{k+1}(({\M}(X)) \ar[d]^{\gamma}
& {\SD}^{k+1}(\SV(X))\ar[d]^{\eta}\ar[l]_{\quad\iota'}^{\quad\sim}\\ 
{\CD}^b(C\M^k(X))\ar[r]_{\zeta}^{\sim} &{\CD}^b({\M}^k(X))\ar[r]_{\iota}^{\sim}
& {\SD}^k(({\M}(X)) 
& {\SD}^k(\SV(X))\ar[l]^{\quad\iota'}_{\quad\sim}\\\ 
\enddiagram 
\end{equation} 
Then, all the horizontal functors are equivalences of derived categories and all the  vertical functors are fully faithful.
\eT
\pf The equivalences of the horizontal functors follows from Theorem \ref{mainThmFoxII} and the results in 
\cite[Theorem 3.2]{FM}. It also follows from Theorem \ref{mkp1toMkFlF} that $\beta$ is fully faithful. This completes the proof.
 $\eop$ 

\section{Implications in ${\bf K}$-theory and others} 
\label{AppSection}

In this section we discuss the implications of the equivalence Theorem 
\ref{equivThm}. 
We will not repeat the prelude we provided in the introduction. 
First, we consider the consequences  in $K$-theory.
Our standard reference for $K$-theory would be \cite{S1} and we freely use the definitions and notations from \cite{S1}.
 However, for an exact category $\SE$ or a complicial exact category $\SE$ with weak equivalences,
  ${\K}(\SE)$ will denote the $\K$-theory spectra of $\SE$ and $\K_i(\SE)$ will denote the $\K$-groups. 
Likewise, ${\bf K}(\SE)$ would denote the ${\bf K}$-theory space of $\SE$.
First, we recall a notation and a lemma.  For a noetherian scheme $X$, denote
  $$
  X^{(k)}:=\left\{Y\in X: co\dim\left(Y\right)=k \right\}
  \quad {\rm and~recall}\quad X_x:=\spec{{\CO}_{X,x}}.
 $$
 We recall the following well known result that follows from \cite{Bs, BW}.
 
 \bL\label{uptoFactor}
  Suppose $X$ is a noetherian quasi-projective scheme over an affine scheme $\spec{A}$ and $k\geq 0$ is an integer.
Then, the sequence of derived categories
$$
\diagram 
{\SD}^{k+1}\left(\SV(X)\right) \ar[r] &  {\SD}^{k}\left(\SV(X)\right) \ar[r] & \coprod_{x\in X^{(k)}}{\CD}^{k}\left(\SV(X_x)\right) \\
\enddiagram
$$
 is exact up to factor. If $X$ is regular, this sequence is exact. 
 
 \eL

 \subsection{$K$-theory}\label{SecKTheo}
 
 The following is the main application of (\ref{equivThm}) to ${\K}$-theory. 
 \bT
 \label{Ktheorem}
 Suppose $X$ is a noetherian quasi-projective scheme over an affine scheme $\spec{A}$ and $k\geq 0$ is an integer.
 Consider the diagram of $\K$-theory spectra and maps:
$$
\diagram 
\K\left(C\M^{k+1}(X)\right) \ar[d]_{\wr} & \K\left(C\M^k(X)\right) \ar[d]_{\wr}^{\Psi} &\coprod_{x\in X^{(k)}} \K\left(C\M^k\left(X_x \right)\right)\ar[d]^{\wr}\\
\K\left({\SC}h^{k+1}\left(\M(X)\right), \SQ\right) \ar[r] & \K\left({\SC}h^{k}\left(\M(X)\right), \SQ \right) \ar[r] &\coprod_{x\in X^{(k)}} \K\left({\SC}h^{k}\left(\M\left(X_x\right)\right), \SQ  \right)\\
\K\left({\SC}h^{k+1}\left(\SV(X)\right), \SQ \right) \ar[r] \ar[u]^{\wr}& \K\left({\SC}h^k\left(\SV(X)\right), \SQ \right) \ar[r] \ar[u]^{\wr}_{\Phi}
&\coprod_{x\in X^{(k)}} \K\left({\SC}h^k\left({\SV}\left(X_x\right)\right), 
\SQ  \right)\ar[u]_{\wr}\\
\enddiagram 
$$
Then,the vertical maps are  homotopy equivalences of 
 ${\K}$-theory spectra. 
 Further, if $X$ is Cohen-Macaulay, then the second line and the third line are  homotopy fibrations of ${\K}$-theory spectra. 
\eT
\pf Here the middle upward arrow $\Phi:\K\left({\SC}h^k\left(\SV(X), \SQ \right)\right)
\to \K\left({\SC}h^{k}\left(\M(X), \SQ \right)\right)$
is induced by the functor 
$\iota':\left({\SC}h^k\left(\SV(X), \SQ \right)\right)
\to \left({\SC}h^{k}\left(\M(X), \SQ \right)\right)$ of complicial exact categories with weak equivalences. 
By (\ref{equivThm}), $\iota'$ induces an equivalence of the associated triangulated categories. Therefore, 
by (\cite[3.2.29]{S1}) $\Phi$ is a homotopy equivalence. Likewise, other two upward arrows are homotopy equivalences. 

The middle downward arrow $\Psi$ 
is a composition of three maps,  as follows:
$$
\diagram 
\K\left(C\M^k(X)\right) \ar[d]_{\Psi} \ar[d]\ar[r]^{\Psi'\qquad}&\K\left(Ch^b\left(C\M^k(X), \SQ \right)\right)\ar[d]^{\zeta'}\\
\K\left({\SC}h^{k}\left(\M(X), \SQ \right)\right) &\K\left(Ch^b\left(\M^k(X), \SQ \right)\right) \ar[l]^{\iota'}\\
\enddiagram 
$$
Now, $\zeta'$ and $\iota'$ are induced by the corresponding functors of complicial exact categories, with weak equivalences. Again,
by (\cite[3.2.29]{S1}) in conjunction with Theorem \ref{equivThm}, $\zeta'$ and $\iota'$ are homotopy equivalences. Now, $\Psi'$
is a homotopy equivalence by the agreement theorem \cite[3.2.30]{S1}. Hence, so is $\Psi$.

It remains to show that, when $X$ is Cohen-Macaulay,  the third line is a homotopy fibrations of $\K$-theory spectra. To do this, 
consider sequence of complicial exact categories with weak equivalences (not necessarily exact): 
%
\begin{equation}\label{exatCE2}
\diagram 
\left({\SC}h^{k+1}\left(\SV(X)\right), \SQ \right) \ar[r] & \left({\SC}h^{k}\left(\SV(X)\right), \SQ \right)  \ar[r] &\coprod_{x\in X^{(k)}}
\left({\SC}h^{k}\left(\SV(X_x)\right), \SQ \right) \\
\enddiagram 
\end{equation} 
By (\ref{uptoFactor}), the corresponding sequence of the derived categories is exact up to factor.
Therefore, by an application of the non-connective version of the 
Thomason-Waldhausen localization theorem (see \cite[3.2.27]{S1}) the 
third line in the statement of the theorem is 
a homotopy fibration of $\K$-theory spectra. 
\pic $\eop$  

The following is an immediate consequence of (\ref{Ktheorem}).

\bC\label{BKiCM}
Let $X$ and $k$ be as in Theorem \ref{Ktheorem}.
Assume $X$ is Cohen-Macaulay. 
Then, for any integer $n$, 
there is an exact sequence of $\K$-groups,
$$
\diagram
\cdots \ar[r] &\K_n\left(C\M^{k+1}(X)\right) \ar[r] & \K_n\left(C\M^{k}(X)\right)\ar[r] &  \oplus_{x\in X^{(k)}}\K_n\left(\M^{k}(X_x)\right)\\
\ar[r] & \K_{n-1}\left(C\M^{k+1}(X)\right) \ar[r] & \cdots&\\
\enddiagram
$$  
\eC
\pf Follows from Theorem \ref{Ktheorem}.
\pic $\eop$

\begin{remark}\label{bfkThm}{\rm If $X$ is regular,
statements exactly similar to Theorem \ref{Ktheorem} and Corollary \ref{BKiCM}, respectively, for ${\bf K}$-theory spaces and groups, would be valid. 
Regularity was used to apply the connective version of Thomason-Waldhausen Localization theorem (\cite[3.2.23]{S1}), which requires that the corresponding sequence of 
derived categories is exact.  
}
\end{remark}

\begin{remark} \label{remGerSPSeq}{\rm 
Let $X$ be as in Theorem \ref{Ktheorem}.
Assume $X$ is Cohen-Macaulay. 
The following are some  remarks. 
\bE
\item The usual diagram to compute the Gersten complex, reduces to
$$
\diagram
\oplus_{x\in X^{(k-1)}}\K_{n+1}\left(C\M^{k-1}(X_x)\right)\ar[d]\ar@{-->}[dr]&&
\K_{n-1}\left(C\M^{k+2}(X)\right)\ar[d]\\
\K_n\left(C\M^{k}(X)\right)\ar[r]\ar[d] &  \oplus_{x\in X^{(k)}}\K_n\left(C\M^{k}(X_x)\right)\ar[r]\ar@{-->}[dr]
 & 
 \K_{n-1}\left(C\M^{k+1}(X)\right) \ar[d]\\
\K_n\left(C\M^{k-1}(X)\right)&&
\oplus_{x\in X^{(k+1)}}\K_{n-1}\left(C\M^{k+1}(X_x)\right)\\ 
\enddiagram
$$ 
The dotted diagonal arrows form the Gerstner complexes. This provides further insight regarding the Gersten complexes in terms of the groups $\K_n\left(C\M^{k}(X)\right)$. 
The complex is analogous to the $G$-theoretic Gerstner complex in (\cite[Proposition 5.8]{Q}). These are clearly non-isomorphic. 
\item The 
spectral sequence given in \cite{Bn} takes the following form:
$$
E_1^{p,q}=\bigoplus_{x\in X^{(p)}}\K_{-p-q}(C\M^p(X_x))\Lra \K_{-n}(\SV(X)) \qquad{\rm along}\quad p+q=n. 
$$ 
\eE
}
\end{remark} 

\subsection{Grothendieck-Witt theory}
\label{secGrotWitt}
In this section, we develop a counter part of the results on $K$-theory (\ref{Ktheorem}), for Grothendieck-Witt theory. 
First, incorporating duality to the Theorem \ref{equivThm}, we obtain the following.

\bP\label{WittFunt}
Let $X$ be a noetherian quasi-projective scheme, over an affine scheme
$\spec{A}$,  and $k\geq 0$ be an integer.  Then, there is a duality preserving equivalence\\
  ${\CD}^b\left(C\M^k(X)\right) \lra T^k\SD^k\left((\SV(X)\right)$ of the derived categories, where the duality on \\
   ${\CD}^b\left(C\M^k(X)\right)$ is induced by ${\CE}xt^k\left(-, \CO_X\right)$ and that on 
   $T^k\SD^k\left((\SV(X)\right)$ is
  $\#:=T^k{\CH}om\left(-, \CO_X\right)$.
\eP 
\pf It is a standard fact that there is a functor $\M(X) \lra {\CD}^b\left(\SV(X)\right)$, by resolution (e.g. see \cite[3.3]{MNoe}).
The restriction to this functor to $C\M^k(X)$ extends to a functor ${\CD}^b\left(C\M^k(X)\right) \to \SD^k\left((\SV(X)\right)$.
It turns out that this functor represents the composite functor in (\ref{equivThm}). Hence the functor is an equivalence. The routine 
checking establishes that this functor preserves the duality, as required.
\pic $\eop$

For clarity, we point out technical differences in the literature between the basis of $K$-theory, derived Witt theory  and Grothendieck Witt theory.
Recall that $K$-theory is available for complicial exact categories with weak equivalences \cite{S1, TT, W} and derived Witt theory was defined for triangulated categories
with duality \cite{TWG1}. However, in \cite{S3}, the Grothendieck Witt theory ($GW$)  were developed for dg categories with weak equivalences and duality. 
Further, note that $GW$-theory of dg categories encompasses Witt theory \cite[Proposition 6.3]{S3}.
As was pointed out in the introduction, $K$-theory is invariant of equivalences of the associated triangulated categories, and when $2$ is invertible,
so are $GW$-theory and Witt theory. 
The primary reason why we cannot directly imitate the methods   in  $K$-theory (\S \ref{SecKTheo}), for Grothendieck Witt theory is that 
the dg category ${\SC}h^k(\M(X))$ 
does not have a natural duality structure. 
The remedy for this is obtained by embedding ${\SC}h^k(\M(X))$ in the 
 dg category of Prefect complexes. 
Following are some notations and background.

\begin{notations}\label{notaPerfect}
{\rm 
We establish some background as follows.
Let $X$ denote a noetherian scheme. 
\bE
\item We will donate the category of quasi-coherent $\CO_X$-modules by $QCoh(X)$. 
Also, $Ch(QCoh(X))$ will denoted the category of chain complexes of objects in $QCoh(X)$
and $\CD(QCoh(X))$ will donated its derived category (see \cite[A.3.2]{S1}). 
\item Recall, that a complex $\CF_{\bullet}\in Ch(QCoh(X))$ is called perfect complex, if for all $x\in X$, there is an affine open
 neighborhood $U$ and a quasi-isomorphism $\CE_{\bullet}\to (\CF_{\bullet})_{|U}$, for some $\CE_{\bullet}\in Ch^b(\SV(U))$. 
 This is equivalent to saying $\CE_{\bullet}$ isomorphic to $(\CF_{\bullet})_{|U}$ in the derived category $\CD(QCoh(X))$ (see \cite[Lemma 2.2.9]{TT}).
  
\item Denote the category of perfect complexes of $\CO_X$-modules by $Perf(X)$. For a closed subscheme $Z$, 
$Perf_Z(X)$ will denote the subcategory of complexes ${\CG}_{\bullet} \in Perf(X)$, such that   
$Supp\left({\CH}_i(\CG_{\bullet})\right)\subseteq Z$. The corresponding derived categories would be,
respectively, denoted by 
${\CD}(Perf(X))$ and ${\CD}(Perf_Z(X))$. 
\item In analogy to notations \ref{nota}(\ref{SCHSDk}), for integers $k\geq 0$, 
$$
\left\{
\begin{array}{l}
{\SP}erf^k(X):=\left\{ 
\CF_{\bullet}\in Perf(X): \forall ~i~\CH_i\left(\CF_{\bullet}\right) \in Coh^k(X)
\right\}\\
{\SD}^kPerf(X):=\left\{ 
\CF_{\bullet}\in {\CD}^b(Perf(X)): \forall ~i~\CH_i\left(\CF_{\bullet}\right) \in Coh^k(X)
\right\}\\
\end{array}
\right.
$$ 
 would denote the full subcategory, of the 
 respective categories, of such objects. Note, ${\SD}^kPerf(X)$ is the derived category of 
 ${\SP}erf^k(X)$. 
 \item To avoid confusion, we will use prefix $dg$ to denote the respective dg categories. So, 
 $dgPerf(X)$ would denote the dg category whose objects are same as that of $Perf(X)$. Likewise,
 $dg\SV(X)$, $dg^k\SV(X)$, $dgCM^k(X)$ will, respectively, denote the dg categories 
 whose objects are, respectively, the same as that of  $Ch^b(\SV(X))$, ${\SC}h^k(\SV(X))$, $Ch^b(C\M^k(X))$.
 \item \label{itemIDot} Throughout, we fix a minimal injective resolution $I_{\bullet}$ of $\CO_X$, as follows:
 $$
 \diagram 
 0\ar[r] & \CO_X \ar[r] & I_0\ar[r] & I_{-1} \ar[r] & I_{-2}\ar[r] &\cdots\\
 \enddiagram.
 \qquad {\rm Clearly,}~ I_{\bullet}\in Perf(X). 
 $$
 For 
$\CF_{\bullet} \in Perf(X)$, denote $\CF^{\vee}:={\CH}om(\CF_{\bullet}, I_{\bullet})$. 
For properties of such  minimal resolutions and the nature of arguments, the readers are referred to \cite{BH} and \cite{G}. 
\eE 
}
\end{notations} 
 The following addresses the duality aspect of $dgPerf(X)$.
 
 \bL\label{dgDuality}
 Let $X$ be a noetherian scheme. Let $I_{\bullet}$ be as in Notation \ref{notaPerfect} (\ref{itemIDot}).
 Then, the association $\CF_{\bullet} \mapsto \CF_{\bullet}^{\vee}$ endows  $(dgPerf(X), \SQ)$ with a structure of 
 a dg category with weak equivalences and duality, weak equivalences being  the set of all quasi-isomorphism $\SQ$.
 \eL
 \pf Consider $\CO_X$ as a complex, concentrated at degree zero. Since $\CO_X\to I_{\bullet}$ is a quasi-isomorphism,
 $I_{\bullet}\in Perf(X)$ . Let $\CF_{\bullet} \in Perf(X)$. So, there is an affine open subset $U$ and a
 quasi-isomorphism  $\CE_{\bullet}\to (\CF_{\bullet})_{|U}$ for some $\CE_{\bullet}\in Ch^b(\SV(U))$.
 Then,
 $$
 {\rm In}~\CD(QCoh(X)),\quad
 \diagram
 {\CH}om(\CE_{\bullet}, \CO_U) \ar[r]_{\sim~~~} & {\CH}om(\CE_{\bullet}, (I_{\bullet})_{|U})& {\CH}om((\CF_{\bullet})_{|U}, (I_{\bullet})_{|U})\ar[l]_{\sim~~} \\
 \enddiagram
 $$
 Since ${\CH}om(\CE_{\bullet}, \CO_U) \in Ch^b(\SV(U))$, it follows, $\CF_{\bullet}^{\vee}\in Perf(X)$.
 Also, for $\CF_{\bullet}\in Perf(X)$, we check, $\CF_{\bullet} \to \CF_{\bullet}^{\vee\vee}$ is a quasi-isomorphism. 
 As above, let  $U$ be affine and $\CE_{\bullet}\to (\CF_{\bullet})_{|U}$ be a quasi-isomorphism, where $\CE_{\bullet}
 \in Ch^b(\SV(X))$. Since, the question is local, we can assume that $X=U$.
 Consider the diagram
 $$
 \diagram
 \CE_{\bullet}\ar[d] \ar[rr] && {\CH}om\left( {\CH}om\left(\CE , \CO_U\right), \CO_U\right)\ar[d]\\
 \CF_{\bullet} \ar[rr] && \CF_{\bullet}^{\vee\vee}\\
 \enddiagram
 $$
Since, two vertical arrows and the top horizontal arrows are quasi-isomorphisms, so is the bottom horizontal arrow. 
Now, we show if $f_{\bullet}:\CF_{\bullet}\to \CG_{\bullet}$ is a quasi-isomorphism in $Perf(X)$, so is
 $f_{\bullet}^{\vee}:\CG_{\bullet}^{\vee}\to \CF_{\bullet}^{\vee}$. Using similar arguments as above, we can 
 assume that $X$ is affine and $f$ is map in $Ch^b(\SV(X))$, in which case the assertion is well known.
 \pic  $\eop$
 
The following proposition on derived equivalences is derived from results in \cite{TT}. 
\bP\label{ttDerEq}
Let $X$ be noetherian separated scheme, with an ample family of line bundles and $k\geq 0$ be an integer.
Let $I_{\bullet}$ be as in Notation \ref{notaPerfect} (\ref{itemIDot}).
Then, 
$\SD^k(\SV(X))\to \SD^kPerf(X)$ is an equivalence of derived categories. 
\eP
 \pf We have $\CD^b(\SV(X)) \to \CD(Perf(X))$ is an equivalence of derived categories (\cite[Lemma 3.8]{TT}, \cite[Prop. 3.4.8]{S1}).
 Since $\SD^k(\SV(X)) \hra \CD^b(\SV(X))$ and $\SD^kPerf(X) \hra \CD(Perf(X))$ are full subcategories, the assertion follows. 
 $\eop$  
 
 Of our particular interest would be the following equivalences of derived categories.
 \bP\label{ourDGFs}
 Suppose $X$ is a quasi-projective scheme over an affine scheme $\spec{A}$ and $k\geq 0$ is an integer.
 Then, 
 \bE
 \item \label{perfSV}The inclusion functor $dg^k(\SV(X) \hra dg{\CP}erf^k(X)$ is a duality preserving form functor
 {\rm (see \cite[1.12, 1.7]{S3}, for definition)}, of pointed dg categories with weak equivalences and dualities, such that the associated 
 functor of the triangulated categories 
 $\ST\left(dg^k\SV(X)\right) \hra \ST\left(dg{\CP}erf^k(X)\right)$ is an equivalence. 
 \item \label{dpfCM}The inclusion functor $dgC\M^k(X) \hra T^k\left(dgPerf^k(X)\right)$ is a duality preserving form functor, 
 of pointed dg categories with weak equivalences and dualities,
 where $T$ denotes the shift.
Further, the associated 
 functor of the triangulated categories \\
 $\ST(dgC\M^k(X))\to \ST\left(T^k\left(dgPerf^k(X)\right)\right)$ is an equivalence. 
 \eE
 \eP
\pf Since $\ST\left(dg^k(\SV(X)\right)=\SD^k(\SV(X))$ and $\ST\left(dg{\CP}erf^k(X)\right)=\SD^k(Perf(X))$,
the latter part of (\ref{perfSV}) follows immediately from  Proposition \ref{ttDerEq}. The duality compatibility
transformation 
 is the obvious map ${\CH}om(\CF_{\bullet}, \CO_X)  \to {\CH}om(\CF_{\bullet}, I_{\bullet})$,
which is a weak equivalence. This establishes (\ref{perfSV}).

To prove (\ref{dpfCM}), note that the duality on $dgC\M^k(X)$ is induced by the duality ${\CE}xt^k\left(-, \CO_X\right)$. 
For $\CF_{\bullet}\in dg(C\M^k(X))$, let $\widehat{\CF_{\bullet}}$ denote its dual.
Note, $\ST\left(dgC\M^k(X)\right) =\CD^b\left(C\M^k(X)\right)$
and $\ST\left(dg{\CP}erf^k(X)\right)=T^k\SD^k\left(Perf^k(X)\right)$. 
Now, it follows from (\ref{perfSV}) and Proposition \ref{WittFunt}, that 
$\ST(dgC\M^k(X))\to \ST\left(T^k\left(dgPerf^k(X)\right)\right)$ is an equivalence. 
For an object $\CF\in C\M^k(X)$, due to grade consideration it follows that 
$\forall ~j\leq k-1, ~{\CH}om(\CF, I_j)=0$ (see \cite[Prop. 3.2.9]{BH}, \cite[Theorem 1.15]{G}). Therefore, for a complex $\CF_{\bullet}\in 
dgC\M^k(X)$, we have a bounded double complex:
$$
\diagram
& 0 \ar[d] &0 \ar[d] & 0\ar[d] & \\
\cdots \ar[r] & {\CE}xt^k\left(\CF_{n+1}, \CO_X\right) \ar[r] \ar[d]& {\CE}xt^k\left(\CF_{n}, \CO_X\right) \ar[r] \ar[d]& {\CE}xt^k\left(\CF_{n-1}, \CO_X\right) \ar[r] \ar[d]& \cdots\\
\cdots \ar[r] & {\CH}om\left(\CF_{n+1}, I_{-k}\right) \ar[r] \ar[d]&  {\CH}om\left(\CF_{n}, I_{-k}\right) \ar[r] \ar[d]& {\CH}om^k\left(\CF_{n-1}, I_{-k}\right) \ar[r] \ar[d]& \cdots\\
\cdots \ar[r] & {\CH}om\left(\CF_{n+1}, I_{-(k+1)}\right) \ar[r] \ar[d]&  {\CH}om\left(\CF_{n}, I_{-(k+1)}\right) \ar[r]\ar[d] & {\CH}om^k\left(\CF_{n-1}, I_{-(k+1)}\right) \ar[r] \ar[d]& \cdots\\
& \cdots &\cdots & \cdots & \\
\enddiagram
$$
In this double complex, the vertical lines are exact. This gives a natural transformation 
$\varphi: \widehat{\CF_{\bullet}}\to {\CH}om\left(\CF_{\bullet}, I_{\bullet}\right)=:\CF_{\bullet}^{\vee}$, which is a weak equivalence.
Now, it follows that $\varphi$ defines a duality preserving transformation $dgC\M^k(X)\to T^k\left(dgPerf^k(X)\right)$.
\pic $\eop$ 

The following useful diagram is analogous to the diagram in the Equivalence Theorem \ref{equivThm}, in the context of 
dg categories with weak equivalences and dualities. 
\bC\label{dgEuivDia}
Suppose $X$ is a quasi-projective scheme over an affine scheme $\spec{A}$, and $k\geq 0, r$ are integers.  
Consider the diagram
\begin{equation}\label{diaDGEqThm}
\diagram
T^{-1}dgC\M^{k+1}(X)  \ar[r] &T^{k}dg{\CP}erf^{k+1}(X) \ar[d] &T^{k}dg^{k+1}\SV(X) \ar[l]\ar[d]\\
dgC\M^k(X)  \ar[r] &T^kdg{\CP}erf^k(X)   & T^kdg^k\SV(X)\ar[l]\\
\enddiagram 
\end{equation}
In this diagram, all the arrow are form functors of dg categories with weak equivalence. Further, the horizontal arrows induce 
equivalences of the associated triangulated categories and the right hand square commutes. Note that there is no natural 
vertical functor on the left side. 
\eC
\pf Follows from Proposition \ref{ourDGFs}. 
$\eop$

Now we have the machinery to state our results on $GW$-theory. For further background information regarding the definitions of 
${\bf GW}$-spectrum and ${\G}W$-spectrum the readers are referred to Section \ref{secSpGWKGW} oe \cite{S3}. 
\bT\label{GWThm}
Suppose $X$ is a quasi-projective scheme over an affine scheme $\spec{A}$, with $1/2\in A$ and $k\geq 0, r$ are integers.  
In the following, weak equivalences and dualities in the respective categories would be as in (\ref{ourDGFs}).
Then, the maps in the following zig-zag sequences
%
$$
\diagram 
{\bf GW}^{[r]}\left(dgC\M^k(X)\right) \ar[r]^{\zeta}& {\bf GW}^{[k+r]}\left(dgPerf^k(X)\right) &{\bf GW}^{[k+r]}\left(dg^k\SV(X)\right)\ar[l]_{\Phi} & {\rm in}~ \Sp\\
{\G}W^{[r]}\left(dgC\M^k(X)\right) \ar[r]_{\zeta} & {\G}W^{[k+r]}\left(dgPerf^k(X)\right)&{\G}W^{[k+r]}\left(dg^k(\SV(X)\right)\ar[l]^{\Phi} & {\rm in}~ \mathrm{BiSp}\\
\enddiagram 
$$
are stable homotopy equivalences in the respective categories. 
\eT
\pf Follows directly from Proposition \ref{ourDGFs} and  \cite[Theorem 6.5]{S3},  \cite[Theorem 8.9]{S3}. 
\pic $\eop$

\bT\label{HomFibrationGWs}
Suppose $X$, $k, r$  are as in (\ref{GWThm}). Assume further that  $X$ is a Cohen-Macaulay scheme. 
Consider the following diagram of ${\G}W$-spetra:
$$
\diagram 
{\G}W^{[-1+r]}\left(dgC\M^{k+1}(X) \right) \ar[d] &
{\G}W^{[r]}\left(dgC\M^{k}(X)  \right) \ar[d] &
\coprod_{x\in X^{(k)}}{\G}W^{[r]}\left(dgC\M^{k}\left(X_x\right)  \right) \ar[d] \\
{\G}W^{[k+r]}\left(dgPerf^{k+1}(X)  \right)  &
{\G}W^{[k]}\left(dgPerf^{k}X)  \right)  &
\coprod_{x\in X^{(k)}}{\G}W^{[k+r]}\left(dgPerf(X_x)\right)  \\
{\G}W^{[k+r]}\left(dg^{k+1}\SV(X) \right) \ar[r] \ar[u]&
{\G}W^{[k+r]}\left(dg^{k}\SV(X) \right) \ar[r] \ar[u]&
\coprod_{x\in X^{(k)}}{\G}W^{[k+r]}\left(dg^k\SV(X_x) \right) \ar[u] \\
\enddiagram
$$
In this diagram, all the vertical arrows are equivalence of homotopy Bispectra and the bottom sequence is a homotopy fibration of bispectra. 
Further, if $X$ is regular then the corresponding statement for ${\bf GW}$-spectra would be valid. 
\eT
\pf It follows directly from Theorem \ref{GWThm} that the vertical rows are equivalences. It remains to show that the bottom row is 
a fibration. This follows from the localization theorem \cite[Thm 8.10]{S3} and Lemma \ref{uptoFactor}. When $X$ is regular, use 
the localization theorem \cite[Thm 6.6]{S3}. \pic 
$\eop$ 
\begin{remark}{\rm
The following are some remarks:
\bE
\item 
As in Corollary \ref{BKiCM}, each shift $r$, corresponding to the fiber sequence in Theorem \ref{HomFibrationGWs}, an exact sequence of 
${\G}W$-groups would follow. Likewise, analogous to Remark \ref{remGerSPSeq}, for each shift $r$, a spectral sequence ${\G}W$-groups would follow.
\item For a scheme $X$ and a rank one  locally free sheaf $\CL$, ${\CH}om(-, \CL)$ induces a duality on $\SV(X)$. All off the above would be valid, with 
dualities induced by ${\CH}om(-, \CL)$, instead of ${\CH}om(-, \CO_X)$.
\item For an exact category $\SE$ with duality, Schlichting \cite{S4} defined Grothendieck-Witt space $GW(\SE)$. By Agreement theorem 
\cite[Proposition 6]{S4} and \cite[Proposition 5,6]{S3}, with $X$ as in (\ref{GWThm}) and integers
 $k\geq 0$ we have$GW(C\M^k(X))$ is naturally equivalent to the infinite look space $\Omega^{\infty}{\bf GW}(dgC\M^k(X))$.
\eE
}
\end{remark}
\subsection{Derived Witt Theory}
\label{secDerWitt}

%
In this subsection we comment on Witt theory. The following follows immediately.
\bT\label{WCM}
Let $X$, $A$ and $k$ be as in (\ref{WittFunt}). Assume $1/2\in A$.
Then, the maps of the shifted Witt groups
 $W^{r}\left({\CD}^b\left(C\M^k(X)\right)\right) 
\to W^{k+r}\left(\SD^k\left((\SV(X)\right)\right)$ are isomorphisms, for
all $r\in \Z$.
In particular, the maps
$$
\left\{
\begin{array}{l}
W\left(C\M^k(X)\right)\to W^{k+4r}\left(\SD^k\left((\SV(X)\right)\right)\\
W^-\left(C\M^k(X)\right)\to W^{k+2+4r}\left(\SD^k\left((\SV(X)\right)\right)\\
\end{array}
\right.
$$
are isomorphisms, for all $r\in \Z$. 
\eT
\pf The first isomorphism follows from (\ref{WittFunt}) and \cite[Theorem 6.2]{TWG1}, because the quotient category would be trivial.
By \cite[Theorem 1.4]{BW}, $W\left(C\M^k(X)\right)\cong W^{4r}\left({\CD}^b\left(C\M^k(X)\right)\right)$  and
$W^-\left(C\M^k(X)\right)\cong W^{2+4r}\left({\CD}^b\left(C\M^k(X)\right)\right)$. By combining these, with the first isomorphism, 
the latter two isomorphisms are established.
$\eop$

Due to non-availability of Thomason-Waldhausen  \cite[3.2.27]{S1} type of localization theorems in derived Witt theory, 
for the statement of the following theorem, 
 we would assume 
that $X$ is a regular. 

\bT\label{witt5Term}
Let $X$ be a quasi-projective regular scheme, over an affine scheme
$\spec{A}$, with $1/2\in A$, and $k\geq 0$ be an integer. 
In this case, the sequence in Lemma \ref{uptoFactor} is exact.
Now, the twelve term exact sequence in  
{\rm \cite[Corollary 6.6]{TWG1}}, corresponding to 
the same exact sequence  
reduces to two five term exact sequences of Witt groups as follows:
$$
\diagram 
0\ar[r] & W^{-1}\left({\CD}^b\left(C\M^{k+1}(X)\right) \right)  \ar[r] & 
W\left(C\M^{k}(X)\right)  \ar[r] &\oplus_{x\in X^{(k)}} W\left(\M^k(X_x) \right) \\
\ar[r]  & W\left(C\M^{k+1}(X)\right)  \ar[r] &W^{1}\left({\CD}^b\left(C\M^{k}(X)\right) \right) \ar[r]&0\\
\ar[r]&W^1\left(C\M^{k+1}(X)\right)\ar[r]&W^-\left(CM^k(X)\right) \ar[r]&\oplus_{x\in X^{(k)}} W^-\left(\M^k(X_x) \right)\\
\ar[r]&W^-\left(C\M^{k+1}(X)\right)\ar[r]&W^3\left({\CD}^b\left(C\M^b(X)\right) \right)
\ar[r]&0\\
\enddiagram 
$$
\eT
\pf Write down the twelve term exact sequence (\cite[Corollary 6.6]{TWG1}) of Witt groups,
corresponding exact sequence of triangulated categories  above.
The zero term in the second row corresponds to $\oplus_{x\in X^{(k)}}W^{k+1}(X_x)=0$ by \cite{BW}. 
Likewise, the first and the last zero are established. 
The rest follows by identifying the other terms by Corollary \ref{WCM}. $\eop$

%
\appendix
\section{Background  on ${\bf GW}$-Spectrum and ${\G}W$-Bispectrum}\label{secSpGWKGW}
In this section,  we 
include the some background information on ${\bf GW}$-spectrum and ${\G}W$-bispectrum.
 In fact this is an overview from our primary source \cite{S3} on the same.
Readers reluctant to deal with language of model categories may like to refer directly to formulas \ref{coLimPiIGW},  \ref{piColimGW}. 
 
 Recall that for  exact categories $\SE$ and also for complicial exact categories $\SE$ with weak equivalences
the  ${\bf K}$-theory spaces ${\bf K}(\SE)$ 
were 
defined as pointed topological spaces. Further, the   $\K$-theory spectra $\K(\SE)$ were defined as $\Omega$-spectra of topological 
spaces, which is a sequence of pointed topological spaces, with bonding maps (see \cite[\S A.1.8]{S1}). 
 Likewise, for a pointed dg category $\SA$ with weak equivalences and duality
two invariants are defined in $GW$-theory, in \cite{S3}, as follows.
\bD\label{DEFSpGW}
{\rm 
For a dg category $\SA$ with weak equivalences and duality, 
the Grothendieck-Witt spectrum ${\bf GW}(\SA)$ takes value in the the category of symmetric spectra of pointed topological spaces.
We only give some outlines of  the definitions of symmetric spectrum (see \cite[\S B.1]{S3} for details) and  of ${\bf GW}(\SA)$ as follows.
\bE
\item 
Let $\Sigma$ denote the category,  whose objects are $\underline{n}:=\{1, 2, \ldots, n\}$, $\underline{0}=\phi$. The morphisms are:  
$\forall ~m\neq n~Mor(\underline{m}, \underline{n})=\phi$ and 
$\forall~n~Mor(\underline{n}, \underline{n})=:\Sigma_n$ 
is the group of permutations elements in  $\underline{n}$.\\
The category of pointed 
topological spaces is denoted by $\mathrm{Top}_*$. The smash product $\wedge$ provides the category $\mathrm{Top}_*$  a
structure of a symmetric 
monoidal category.  
Also, denote $S^0:=\{0, 1\}$, and $S^{1}$ denotes the circle obtained by identifying $0\sim 1$ in the unit interval $[0,1]$.
For $n\in \N$, denote $S^n=S^1\wedge \cdots \wedge S^1$, the $n$-fold smash product. 
  \item 
A symmetric spectrum is a functor $\Sigma \to Top_*$ form $\Sigma$ to the category $\mathrm{Top}_*$, together with base point preserving maps 
$\forall ~n, m\in \N\quad  e_{n,m}: S^n \wedge X_m \to X_{n+m}$,
to be called the bonding maps,
  %
 with further compatibility conditions. Therefore, a symmetric spectrum is a sequence $X:=\{X_0, X_1, X_2, \ldots\}$ of pointed topological spaces such that 
 (1)~ $\forall~n\in \N$, there is a continuous base point preserving left action of $\Sigma_n$ on $X_n$, and
 (2) for $n, m\in \N$, there are pointed continuous $\Sigma_m\times \Sigma_m$-equivariant maps 
  $
  e_{n,m}: S^n \wedge X_m \to X_{n+m}
  $ with natural compatibility conditions.

  The category of symmetric spectrum is denoted by $\mathrm{Sp}$. The smash product $\wedge$ of pointed topological spaces extends to a smash product on
 $\mathrm{Sp}$, denoted by the same notation $\wedge$. Further, 
 four different  model structures on $\Sp$ is discussed in \cite[\S B.2.]{S3}, namely, the (positive) projective level model structure and 
 (positive) projective stable model structure. Two stable model structures on $\Sp$ have same weak equivalences and hence same homotopy category,
 to be called the stable homotopy category. 
 
 For $X, Y\in \Sp$, let $[X, Y]$ denote the set of all morphisms $X\to Y$ in the stable homotopy category.
 Define the homotopy groups of
 $X$  as 
 $$
 \pi_n(X):=[S^n, X] 
 $$
 Our interest, with respect to $GW$-theory, would remain limited to the case when the spectrum is a (positive) $\Omega$-spectrum, meaning that the 
 bonding maps
 $X_n\to \Omega X_{n+1}$ are weak equivalences of pointed topological spaces. In this case (see \cite[\S B.3.]{S3}), 
 $$
 \forall ~n\in \Z \qquad \pi_n(X)={\mathrm co}\lim_k \pi_{n+k}(X_k).
 $$
 where $\pi_{n+k}(X_k)$ denote the usual homotopy groups.
 \item For a pointed dg category $\SA$, the Grothendieck-Witt spectrum 
  $GW(\SA)$ is defined as a 
  symmetric spectrum (\cite[\S 4.4, Definition 5.4]{S3})
  $$
  {\bf GW}(\SA):=\{GW(\SA)_0, GW(\SA)_1, GW(\SA)_2, \ldots\}.
  $$ 
Further, for $n\in \Z$ the $n-$shifted $GW$-spectrum is defined to be 
$$
{\bf GW}^{[n]}(\SA):={\bf GW}\left(\SA^{[n]}\right)\quad {\rm where}~\SA^{[n]} ~{\rm denote~the}~n{\rm -shifted~dg~category~of}~\SA.
$$
In fact, $ {\bf GW}(\SA)$,  is a positive $\Omega$-spectrum (\cite[Theorem 5.5]{S3}). For $n, i\in Z$, denote  
\begin{equation}\label{coLimPiIGW}
{\bf GW}_i^{[n]}(\SA):=\pi_i\left({\bf GW}^{[n]}(\SA)\right)={\mathrm co}\lim_k \pi_{n+k}\left({\bf GW}^{[n]}(\SA)_k\right),
\end{equation}
where the latter equality is a property of positive $\Omega$-spectra. 
\eE
}
\eD

With respect to sequences $\diagram (\SA_0, w)\ar[r] & (\SA_1, w) \ar[r] & (\SA_2, w) \enddiagram$ of pointed dg categories with weak equivalences and dualities,
  ${\bf GW}$-spectra behave well (see \cite[Theorem 6.6]{S3}), when the associated sequence of triangulated categories are exact (assuming $2$ is invertible).
  However, while dealing with non-regular schemes, as noted in Lemma \ref{uptoFactor}, the relevant sequence of associated 
  triangulated categories are exact only up to factors. To remedy this situation, in analogy to $\K$-theory spectra, Karoubi-Grothendieck-Witt 
  spectrum ${\G}W(\SA)$ of pointed dg categories $\SA$ with weak equivalences and dualities are defined \cite{S3}.   
\bD\label{DEFKGWBiSp}
{\rm 
For dg categories $\SA$ with weak equivalences and dualities,
the Karoubi-Grothendieck-Witt 
  spectrum ${\G}W(\SA)$ takes value in the category of Bispectra. Again, we only outline the definition of 
  Bispectra and ${\G}W(\SA)$ from \cite{S3}. 
\bE
\item Note that the category $(Top_*, \wedge, S^0)$ of pointed topological spaces, with smash product $\wedge$
 is a symmetric monoidal category, where 
the unit is $S^0$. The process used to obtain the category $\mathrm{Sp}$ from $(Top_*, \wedge, S^0)$ 
is fairly formal and is described in \cite[\S ~B.9.]{S3}, where $S^1$ had a special role to play. 
The category of Bispectra is obtained, by iterating the  same process on $\Sp$, as follows.
\bE
\item Denote $S:=\{S^0, S^1, \ldots \}\in \Sp$. For $X, Y\in \Sp$ define a new smash product $X\wedge_SY$ by push forward
$$
\diagram
X\wedge S \wedge Y \ar[r]^{\quad1_X\wedge -} \ar[d]_{-\wedge 1_Y}& X \wedge Y\ar[d] \\
X \wedge Y \ar[r] & X \wedge_S Y\\
\enddiagram
$$
Then, $(\Sp, \wedge_S, S)$ is cofibrantly generated closed symmetric monoidal model category, with the 
positive stable model structure on $\Sp$. By abuse of notations, write $S^1:=S\wedge(S^1, pt, pt, \ldots)
=(S^1, S^1\wedge S^1, S^2\wedge S^1, \cdots)$. Also, let $\tilde{S}^1$ be a cofibrant replacement 
of $S^1\in \Sp$.

\item The process mentioned above \cite[\S B.9.]{S3} is applied to this category $(\Sp, \wedge_S, S)$ with 
the special role played by $\tilde{S}^1\in \Sp$ (see \cite[\S B.11.]{S3}). The category thus obtained is called category $\tilde{S}^1{\rm -}S$-bispecta
which is denoted by 
$$
\mathrm{BiSp}:=\Sp(\Sp, \tilde{S}^1):=\Sp\left((\Sp, \wedge_S, S), \tilde{S}^1\right)
$$
We describe $Bi\Sp$ as follows (see \cite[\S B.11.]{S3}):
\bE
\item Let $\Sp^{\Sigma}$ denote the category of functors $\Sigma \to \Sp$. So, an object in $\Sp^{\Sigma}$ is a sequence 
$(\CX_0, \CX_1, \CX_2, \ldots, )$ where $\CX_n\in \Sp$ are spectra, with a left action of the symmetric groups $\Sigma_n$. 
\item The smash product $\wedge_S$ extends to a smash product $\wedge_S$ in $\Sp^{\Sigma}$. 
\item  Write $\tilde{S}:=(S^0, \tilde{S}^1, \tilde{S}^1\wedge_S \tilde{S}^1, \tilde{S}^1\wedge_S \tilde{S}^1\wedge_S \tilde{S}^1, \ldots)$. 
\item The objects in $ \mathrm{BiSp}$ are $\tilde{S}$-modules $M\in \Sp^{\Sigma}$. This means that there is map (natural transformation)
 $\tilde{S}\wedge M \to M$, compatible with the action of $\Sigma$. 
\item The objects $\CX\in \mathrm{BiSp}$ are also called $\tilde{S}^1{\rm -}S^1$-bispectrum. 
\eE 
\item $\mathrm{BiSp}$ has a (positive) stable symmetric monoidal model structure (see \cite[\S B.9, B.11]{S3}), by results of Hovey \cite{Ho}.
 For bispectra $\CX, \CY\in \mathrm{BiSp}$, let $[\CX, \CY]_{H(\mathrm{BiSp})}$ denote the set of all morphisms $\CX\to \CY$ in this stable homotopy category.
 Define the homotopy groups of
 $\CX$  as 
 $$
 \pi_n(\CX):=[\tilde{S}^n, \CX]_{H(Bi\Sp)} 
 $$
If $\CX=(\CX_0, \CX_1, \CX_2, \ldots)\in \mathrm{BiSp}$ is a level fibrant semistable $\tilde{S}^1{\rm -}S^1$-bispectrum, then (\cite[Lemma B.16]{S3})
\begin{equation}\label{piColimGW}
 \pi_n(\CX):=[\tilde{S}^n, \CX]_{H(\mathrm{BiSp})} =\mathrm{co}\lim_k\left(\pi_{n+k}(\CX_k)\right)
 \end{equation}
\eE 
\item Suppose $\SA$ is a pointed dg category with weak equivalence and duality. Let $\CS^n\SA$ denote the iterated $n$-fold suspension  of $\SA$
(see \cite[\S 8.1]{S3}). Denote ${\G}W(\SA)_n:={\bf GW}(\CS^n\SA)\in \Sp$, the ${\bf GW}$-spectrum of $\CS^n\SA$.  Then, $\Sigma_n$ 
has a left action on ${\G}W(\SA)_n$. The Karoubi-Grothendieck-Witt spectrum ${\G}W(\SA)$ of $\SA$ is defined as the  
$\tilde{S}^1{\rm -}S^1$-bispectrum (see \cite[\S 8.2]{S3})
$$
{\G}W(\SA):=\left({\G}W(\SA)_0, {\G}W(\SA)_1, {\G}W(\SA)_2, \ldots, {\G}W(\SA)_n, \ldots \right)
\in \mathrm{BiSp}.
$$
The Karoubi-Grothendieck-Witt spectrum is defined as a functor
$
{\G}W: dgCatWD_*\to \mathrm{BiSp}
$
from the category of pointed small dg categories with weak equivalences and dualities to the category of Bispectra. 
For $n\in Z$, the Karoubi-Grothendieck-Witt groups are defined as, and are isomorphic to
$$
{\G}W_n(\SA):=\pi_n\left({\G}W(\SA)\right) \cong
{\mathrm co}\lim_k \pi_{n+k}\left({\G}W(\SA)_k\right)
\cong
{\mathrm co}\lim_k \pi_{n+k}\left({\bf GW}(\CS^k\SA)\right)
$$
$$
\cong
{\mathrm co}\lim_k 
\left({\mathrm co}\lim_{m} \left(\pi_{n+k+m}\left({\bf GW}(\CS^k\SA)\right)_m\right)\right)
$$
The latter isomorphisms follow because the formulas \ref{coLimPiIGW},  \ref{piColimGW} would apply. 
For $n, i\in \Z$  define the $n-$shifted Karoubi-Grothendieck Witt spectrum and groups, respectively, as 
$$
{\G}W^{[n]}(\SA):={\G}W\left(\SA^{[n]}\right)\quad {\rm and} \quad
{\G}W_i^{[n]}(\SA):=\pi_i\left({\G}W^{[n]}(\SA)\right).
$$
  \eE 
  }
\eD

\vspace{5mm}
\noindent{\bf Acknowledgement}: {\it The author would like to express his  appreciation to Marco Schlichting, for his sincere
 academic support over a significant period of time, including the author's recent visits to U. of Warwick.}
%



\end{document}